\documentclass{amsart}
\allowdisplaybreaks[4]

\newcommand{\cL}{\mathcal{L}}

\newcommand{\Z}{\mathbf{Z}}
\newcommand{\R}{\mathbf{R}}
\newcommand{\C}{\mathbf{C}}

\newcommand{{\ba}}{\bf a}
\newcommand{\ve}{\varepsilon}

\newcommand{\ga}{\gamma}

\newcommand{\ra}{\rightarrow}
\newcommand{\Om}{\Omega}
\newcommand{\de}{\delta}
\newcommand{\vp}{\varphi}
\newcommand{\na}{\nabla}
\newcommand{\cd}{\cdot}

\newcommand{\De}{\Delta}
\newcommand{\al}{\alpha}

\newcommand{\nn}{\nonumber}
\newcommand{\ze}{\zeta}
\newcommand{\om}{\omega}

\newtheorem{lemma}{Lemma}{\bf}{\it}
\newtheorem{remark}{Remark}{\it}{\rm}

\newtheorem{theorem}{Theorem}

\numberwithin{theorem}{section}
\numberwithin{lemma}{section}
\numberwithin{equation}{section}
\numberwithin{proposition}{section}
\numberwithin{corollary}{section}
\title[Homogenization of a diffusion perturbed by a vector field ]
{On homogenization of a diffusion perturbed by a periodic reflection invariant vector field}

\author{Joseph G. Conlon }

\address{University of Michigan\\ Department of Mathematics\\ Ann Arbor,
  MI 48109-1109}
\email{conlon@umich.edu}


\begin{document}

\maketitle

\begin{abstract}
In this paper the author studies the problem of the homogenization of a diffusion perturbed by a periodic reflection invariant vector field. The vector field is assumed to have fixed direction but varying amplitude. The existence of a homogenized limit is proven and formulas for the effective diffusion constant are given. In dimension $d=1$ the effective diffusion constant is always less than  the constant for the pure diffusion. In $d>1$ this property no longer holds in general.

\end{abstract}

\section{Introduction}
We consider the problem of the homogenization of a diffusion perturbed by a reflection invariant vector field.  The general set up we have in mind is to understand the limit as $\ve \ra 0$ of the solutions $u_\ve$ to an elliptic equation,
\begin{multline} \label{A1}
- \frac 1{2d} \; \De \; u_\ve(x,\om) - 2b_\ve(x,\om)\partial u_\ve (x,\om)/ \partial x_1
\\
  + u_\ve(x,\om) = f(x), \ x = (x_1,...,x_d) \in \R^d, \ \ \om \in \Om.       
\end{multline}
In (\ref{A1}) the function $f : \R^d \ra \R$ is smooth of compact support and $\Om$ is a probability space.  For simplicity we have assumed that the vector field is always in the $x_1$ direction and hence can be described by the scalar function $b_\ve$.  As $\ve \ra 0$ the field becomes rapidly oscillatory and therefore one might expect that $u_\ve(x,\om)$ converges with probability 1 as $\ve \ra 0$ to a homogenized limit $u(x)$ which is the solution to a constant coefficient elliptic equation,
\begin{equation} \label{B1}
-q(b)\frac{\partial^2 u}{\partial x^2_1} - \ \sum^d_{j=2} \ \frac 1{2d} \frac{\partial^2 u}{\partial x^2_j} + u(x) = f(x), \ \ x \in \R^d.
\end{equation}
The effect of the rapidly oscillating vector field $b_\ve$ is contained in the coefficient $q(b)$ in (\ref{B1}).

In order for a limit $u(x)$ satisfying (\ref{B1}) to exist it is necessary to make assumptions concerning the rapidly oscillating field $b_\ve$.  These are primarily that the distribution functions of the variables $b_\ve(x,\cd)$, $x \in \R^d$, are translation and reflection invariant.  To be specific, we assume that there are translation operators $\tau_x : \Om \ra \Om$, $x \in \R^d$, which are measure preserving and satisfy the group properties $\tau_x \tau_y = \tau_{x+y}$, $x,y \in \R^d, \ \tau_0=$identity. Suppose $b:\Om\ra\R$ is a bounded function.  We then set $b_\ve(x,\om) = b(\tau_{x/\ve} \; \om)$, $x\in \R^d, \; \om \in \Om, \; \ve > 0$.  Such a $b_\ve$ has translation invariant distribution functions and is rapidly oscillating as $\ve \ra 0$.  For $b_\ve$ to satisfy reflection invariance we let $R : \R^d \ra \R^d$ be the reflection operator $R(x_1, ..., x_d) = (-x_1,x_2,...,x_d)$, $x = (x_1,...,x_d) \in \R^d$.  We then require $b : \Om \ra \R$ to satisfy the identities,
\begin{equation} \label{C1}
\left< \prod^n_{i=1} b(\tau_{x_i} \; \cd) \right> = (-1)^n \left< \prod^n_{i=1} b(\tau_{Rx_i} \; \cd) \right>, \ x_i \in \R^d, 1\le i \le n, \ n \ge 1,
\end{equation}
where $\left< \cd\right>$ denotes expectation on $\Om$.  Evidently (\ref{C1}) implies that  $\left<b(\cd)\right> = 0$, so the vector field has no net drift.

A concrete example of an $\Om$ and a $b : \Om \ra \R$  which satisfies (\ref{C1}) is given by taking $\Om$ to be a torus, $\Om = \prod^d_{i=1} [0, L_i]$ with periodic boundary conditions and uniform measure.  The operators $\tau_x : \Om \ra \Om$, $x \in \R^d$, are just translation on $\Om$ and reflection invariance of (\ref{C1}) is guaranteed by the condition,
\begin{equation} \label{D1}
b(x_1,x_2,...,x_d) = - b(L_1 - x_1,x_2,...,x_d), \ \ x = (x_1,...,x_d) \in \Om.
\end{equation}
We shall show that for a discrete version of a periodic $\Om$ with $b$ satisfying (\ref{D1}) a homogenized limit exists with $q(b)$ satisfying $0 < q(b) < \infty$.  For $d=1$ one has $q(b) \le q(0) = 1/2$.  For $d > 1$ it is no longer the case that $q(b) \le q(0) =1/2d$ in general although this does hold for $L_1$ sufficiently small.  One might wish to understand this difference between $d=1$ and $d > 1$ by observing that only in $d > 1$ can one construct nontrivial divergence free vector fields.  The homogenized limit of diffusion perturbed by a divergence free vector field necessarily yields an effective diffusion constant which is larger than the constant for the pure diffusion \cite{fp}.

The homogenization problem considered here appears to have only been studied in the case where $\Om$ is an infinite space for which the variables $b(\tau_x \; \cd), \ x \in \R$, are uncorrelated on a scale larger than $O(1)$.  The problem was introduced by Sinai \cite{s} in a discrete setting.  He proved that in dimension $d=1$ a scaling limit of the random walk corresponding to a finite difference approximation to (\ref{A1}) exists with probability 1 in $\Om$.  The limiting process is strongly subdiffusive.  In a subsequent paper Kesten \cite{k} obtained an explicit formula for the distribution of the scaling limit.  For dimension $d \ge 3$ Fisher \cite{f} and Derrida-L\"{u}ck \cite{dl} predicted that a homogenized limit exists as in (\ref{B1}) with $0 < q(b) < \infty$.  This was proved for sufficiently small $b$ by Bricmont-Kupiainen \cite{bk} and Sznitman-Zeitouni \cite{sz} using a very difficult induction argument.  A formal perturbation expansion for $q(b)$ was obtained in \cite{c2002,c2005} where it was shown that each term of the expansion is finite if $d\ge 3$.  One does not expect the series to converge however.  For $d=1,2$ there are individual terms in the perturbation expansion which diverge.

A main difficulty in the homogenization problem (\ref{A1}), (\ref{B1}) is that when $\Om$ is infinite, good a-priori estimates on the solution to (\ref{A1}) do not hold for all configurations of $b(\cd)$.  In contrast such estimates do hold for divergence form equations with zero drift.  The proof of homogenization in these cases is therefore considerably simpler than for the problem (\ref{A1}), (\ref{B1}).  The first proofs of homogenization for divergence form equations were obtained by Kozlov \cite{ko} and Papanicolaou-Varadhan \cite{pv1} in the continuous case.  K\"{u}nneman \cite{ku} proved a corresponding result for the discrete case.  For non-divergence form equations with zero drift the first proofs in the continuous case were given by Papanicolaou-Varadhan \cite{pv2} and Zhikov-Sirazhudinov \cite{zs}.  Lawler \cite{l} and Anshelevich et al \cite{aks} proved homogenization for a discrete version.  See the books of Bolthausen-Sznitman \cite{bs} for an account of the theory in a discrete setting and of Zhikov et al \cite{zko} for the continuous case.

In this paper we shall be concerned with a discrete version of the homogenization problem described by (\ref{A1}), (\ref{B1}), (\ref{C1}).  Thus the probability space $\Om$ is acted upon by translation operators $\tau_x : \Om \ra \Om$ where now $x \in \Z^d$, the integer lattice in $\R^d$, and satisfy the group properties $\tau_x\; \tau_y = \tau_{x+y}$,  
$\tau_0 =$ identity.  For $i=1,...,d$ let ${\bf e}_i \in \Z^d$ be the element with entry 1 in the ith position and 0 in the other positions.  the discrete equation corresponding to (\ref{A1}) is given by
\begin{eqnarray} \label{E1}
u_\ve(x,\om) &-& \sum^d_{i=1} \ \frac 1{2d} \ \left[ u_\ve(x + \ve{\bf e}_i, \om) + u_\ve(x - \ve{\bf e}_i, \om)\right] \\
&-& b(\tau_{x/\ve}) \left[ u_\ve(x + \ve{\bf e}_1, \om) + u_\ve(x - \ve{\bf e}_1, \om)\right] \nn \\
&+& \ve^2 \; u_\ve(x,\om) = \ve^2 f(x), \ \ x \in \Z^d = \ve \; \Z^d, \ \ \om \in \Om. \nn
\end{eqnarray}
We assume that $b : \Om \ra \R$ satisfies $\sup_\om \ |b(\om)| < 1/2d$, in which case (\ref{E1}) is an equation for the expectation value of a function of an asymmetric random walk.  Hence (\ref{E1}) has a unique bounded solution.  We assume that $b$ satisfies the reflection invariant condition (\ref{C1}) (with $x_i \in \Z^d, 1 \le i \le n$, now).  We also assume that $\Om$ is finite, in which case one can see (Lemma 2.4) that $\Om$ is isomorphic to the integer points on a torus and $b$ has the reflection invariance property (\ref{D1}).  In ${\S 2}$ we prove the following theorem (with $\lfloor \cdot \rfloor$ denoting the integer part):

\begin{theorem}  Assume $\Om$ is a finite probability space and the translation operators $\tau_x : \Om \ra \Om$ are ergodic, $x \in \Z^d$.  Then there exists $q(b), \ 0 < q(b) < \infty$ such that with $u_\ve$ the solution to (\ref{E1}) and $u$ the solution to (\ref{B1}),
\[   \lim_{\ve \ra 0} \ \sup_{x\in \R^d, \om \in \Om} \ |u_\ve(\ve  \lfloor x/ \ve \rfloor, \om) - u(x)| = 0.   \]
\end{theorem}

Suppose now that $\Om$ consists of the integer points on the torus $\prod^d_{i=1}[0, L_i] \subset \R^d$ with periodic boundary conditions.  The reflection invariant condition correspond to (\ref{D1}) is given by
\begin{equation} \label{F1}
b(x_1,x_2,...,x_d) = - b(L_1 - 1-x_1,x_2,...,x_d), \ \ x = (x_1,x_2, ...,x_d) \in \Om.
\end{equation}
For $b$ satisfying (\ref{F1}) we prove in $\S$2, $\S$3 the following results concerning the coefficient $q(b)$ of the homogenized equation (\ref{B1}):

\begin{theorem}  
\begin{itemize}
\item [(a)]  For $d=1$ one has $q(b) \le 1/2$.
\item[(b)] If $d \ge 1$ and $L_1 = 2$ one has $q(b) \le 1/2d$.
\item[(c)] If $d=2$ and $L_1=4$ one has $q(b) \le 1/4$.
\item[(d)] If $d \ge 2$ and $L_1 \ge 6$ is even then there exists $b$ with $q(b) > 1/2d$.
\end{itemize}
\end{theorem}

The proofs of (a), (b), (c) are given in $\S 3$ and are based on applications of the Schwarz inequality.  The proof of (c) is quite lengthy and depends crucially on actual numerical values for a Green's function associated with standard random walk on the integers.  The proof of (d) is given in $\S 2$.  One observes that perturbation theory yields $q(b) = 1/2d + O(|b|^2)$ and that the term $O(|b|^2)$ can be positive.  

In the proof of Theorem 1.2 we use a representation for $q(b)$ in terms of invariant measures for random walk on $\Om$ with drift $b$.  Let $\Om_{d-1}$ consist of the integer points on the $d-1$ dimensional torus $\prod^d_{i=2} [0, L_i] \subset \R^{d-1}$ with periodic boundary conditions.  Setting $L_1 = 2L$ with $L$ an integer we define $\hat \Om$ by
\[    \hat \Om = \left\{ (n,y) : 1 \le n \le L, \ y \in \Om_{d-1} \right\},	\]
whence $\Om$ is the double of $\hat \Om$.  Observe that the boundary of $\partial \hat \Om$ is given by
\[	\partial \hat \Om = \left\{ (1,y), \ (L,y) : y \in \Om_{d-1} \right\}.	\]
Let $\vp^*$ be the invariant measure for random walk on $\hat \Om$ with drift $b(\cd)$ in the ${\bf e}_1$ direction and reflecting boundary conditions on $\partial \hat \Om$.  We define $\psi : \Om_{d-1} \ra \R$ by
\[	\psi(y) = [1/2d - b(1,y)] \ \vp^*(1,y), \ y \in \Om_{d-1}.	\]
Then $q(b)$ is given by the formula,
\begin{equation} \label{G1}
q(b) = \left< \psi \left[ -\De_{d-1} + 4\right]^{-1} \psi_R \right>_{\Om_{d-1}},	
\end{equation}
where $\psi_R$ is defined exactly as $\psi$ but with $b$ replaced by $-b$.  In (\ref{G1}) the expectation $\left< \cd \right>_{\Om_{d-1}}$ is the uniform measure on $\Om_{d-1}$ and $\De_{d-1}$ is the $d-1$ dimensional finite difference Laplacian on functions with domain $\Om_{d-1}$.  The normalization of $\psi$ is chosen so that $q(0) = 1/2d$.  The general formula (\ref{G1}) is proven in $\S 4$.

\section{Proof of Theorem 1.1}
We follow the method introduced in \cite{cn} to obtain homogenized limits.  Thus in (\ref{E1}) we put $u_\ve(x,\om) = v_\ve(x, \tau_{x/\ve} \; \om)$ whence (\ref{E1}) becomes
\begin{eqnarray} \label{A2}
v_\ve(x,\om) &-& \sum^{2d}_{i=1} \frac 1{2d} \left[ v_\ve(x + \ve{\bf e}_i , \tau_{{\bf e}_i }\; \om) + v_\ve(x - \ve{\bf e}_i , \tau_{{-\bf e}_i }\; \om) \right] \\
&-& b(\om) \left[ v_\ve(x + \ve{\bf e}_1 , \tau_{{\bf e}_1 }\; \om) - v_\ve(x - \ve{\bf e}_1 , \tau_{{-\bf e}_1 }\; \om) \right] \nn \\
&+& \ve^2 \ v_\ve(x, \om) = \ve^2  \ f(x), \ \ x \in \Z^d_\ve=\ve\Z^d, \ \ \ \om \in \Om. \nn
\end{eqnarray}
Next we wish to take the Fourier transform of (\ref{A2}).  To show that this is legitimate we first show that the solution  $u_\ve(x,\om)$ of (\ref{E1}) decreases exponentially as $x \ra \infty$.

\begin{lemma}  Suppose $f : \Z^d_\ve \ra \R$ has finite support in the set $\{ x = (x_1,...,x_d) \in \Z^d_\ve : |x| < R \}.$  Let  $u_\ve(x,\om)$ be a bounded solution to $(\ref{A1})$.  Then there are constants $C,K(\ve) > 0$ such that
\begin{equation} \label{B2}
|u_\ve(x,\om)| \le C \exp [K(\ve)(R - |x|) ] \|f\|_\infty, \ \ x \in \Z^d_\ve.
\end{equation}
\end{lemma}
\begin{proof} We write $u_\ve(x,\om) = e^{-kx_1} \; u_{\ve,k}(x,\om)$.  Then from (\ref{E1}) the function $u_{\ve,k}$ satisfies
\begin{eqnarray} \label{E2}
&\ & \frac 1{2d} \sum^{d}_{i=2} \left[ 2u_{\ve,k}(x,\om) - u_{\ve,k}(x + \ve{\bf e}_i ,\om) -  u_{\ve,k}(x - \ve{\bf e}_i , \om)\right]/\ve^2 \\
&+& e^{-k\ve} \left[ \frac 1{2d} + b(\tau_{x/\ve}\; \om)\right] \left\{ u_{\ve,k}(x,\om) -  u_{\ve,k}(x + \ve{\bf e}_1 , \om)\right\}\big/\ve^2 \nn \\
&+& e^{k\ve} \left[ \frac 1{2d} - b(\tau_{x/\ve}\; \om)\right] \left\{ u_{\ve,k}(x,\om) -  u_{\ve,k}(x - \ve{\bf e}_1 , \om)\right\}\big/\ve^2 \nn \\
&+& \left\{ 1 - [\cosh k\ve - 1]/d\ve^2 + 2b(\tau_{x/\ve}\;\om) \sinh \; k\ve/\ve^2\right\} u_{\ve,k}(x,\om) \nn \\
&=& e^{kx_1} \ f(x), \ \ \ \ x \in \Z^d_\ve. \nn
\end{eqnarray}
We may assume wlog that $f$ is nonnegative, whence $u_{\ve,k}$ is also nonnegative.  Suppose $u_{\ve,k}$ attains its maximum at a point $\bar x \in \Z^d_\ve$.  Then we have that
\begin{equation} \label{C2}
\left\{ 1 - [\cosh k\ve - 1]/d\ve^2 + 2b(\tau_{x/\ve} \;\om)\sinh \; k\ve/\ve^2\right\} u_{\ve,k}(\bar x,\om)
\end{equation}
\[	\le \exp \left[ k(\bar x \cd {\bf e}_1) \right] \|f\|_\infty, \ \ |\bar x| < R,	\]
whence it follows that
\begin{equation} \label{D2}
u_\ve(x,\om) \le C \exp \left[ k(\bar x - x) \cd {\bf e}_1 \right] \|f\|_\infty, \ x \in \Z^d_\ve.
\end{equation}
We need to show that  the point $\bar x$ exists for sufficiently small $k$.  To see this assume for contradiction that it does not exist.  Then (\ref{E2}) implies that $\sup_{|x| \le N} u_{\ve,k}(x,\om)$ grows exponentially in $N$ as $N \ra \infty$.  The rate of exponential growth remains bounded away from 0 as $k \ra 0$.  Hence, taking $k$ sufficiently small, we conclude that the function $u_\ve$ is unbounded, contradicting our assumption on $u_\ve$.  The inequality  (\ref{B2}) now follows from (\ref{C2}), (\ref{D2}) on generalizing to all directions ${\bf e}_j, 1 \le j \le d$. 
\end{proof}
For $\xi \in [-\pi/\ve, \; \pi/\ve]^d$ we put
\begin{equation} \label{F2}
\hat v_\ve(\xi, \om) = \int_{\Z^d_\ve} v_\ve(x,\om)e^{ix\cd \xi} dx = \sum_{x\in \Z^d_\ve}\ve^d v_\ve(x,\om)e^{ix\cd\xi}.
\end{equation}
Then from (\ref{A2}) we have that
\begin{eqnarray} \label{G2}
\hat v_\ve(\xi, \om)  &-& \sum^d_{j=1} \frac {1}{ 2d} \left[ e^{-i\ve\xi_j} \hat v_\ve(\xi, \tau_{{\bf e}_j} \om) + e^{i\ve\xi_j} \hat v_\ve(\xi, \tau_{-{\bf e}_j} \om)\right] \\
&-& b(\om)  \left[ e^{-i\ve\xi_1} \hat v_\ve(\xi, \tau_{{\bf e}_1} \om) - e^{i\ve\xi_1} \hat v_\ve(\xi, \tau_{{-\bf e}_1} \om)\right] + \ve^2 \; \hat v_\ve(\xi, \om) \nn \\
&=& \ve^2 \hat f_\ve(\xi), \ \xi \in [-\pi/\ve, \; \pi/\ve]^d, \ \ \om \in \Om, \nn
\end{eqnarray} 
where $\hat f_\ve$ denotes the discrete Fourier transform (\ref{F2}) of $f$.  To solve (\ref{G2}) we define for $\zeta \in [-\pi, \pi]^d$ an operator $\mathcal{L}_\zeta$ on functions $\Psi : \Om \ra \C$ defined by
\begin{eqnarray} \label{H2}
\mathcal{L}_\zeta \Psi(\om) &=& \Psi(\om) - \sum^d_{j=1} \frac 1{2d} \left[ e^{-i\zeta_j} \Psi(\tau_{{\bf e}_j}\om) + 
e^{i\zeta_j} \Psi(\tau_{-{\bf e}_j}\om) \right] \\
&-& b(\om) \left[ e^{-i\zeta_1} \Psi(\tau_{{\bf e}_1}\om) - e^{i\zeta_1} \Psi(\tau_{-{\bf e}_1}\om) \right]. \nn
\end{eqnarray}
Next we define an operator $T_{\eta,\zeta}, \ \eta > 0$,  $\zeta \in [-\pi, \pi]^d$ on $L^\infty(\Om)$ by
\begin{equation} \label{I2}
T_{\eta,\zeta} \; \vp(\om) = \eta \left[ \mathcal{L}_\zeta + \eta\right]^{-1} \vp(\om), \ \ \om \in \Om.
\end{equation}
It is easy to see that $T_{\eta,\zeta}$ is a bounded operator on $L^\infty(\Om)$ with norm at most 1. In fact the RHS of (\ref{I2}) is the expectation for a continuous time random walk on $\Om \times \Z^d$.  The walk is defined as follows:
\begin{itemize}
\item[(a)]  The waiting time at $(\om, x) \in \Om \times \Z^d$ is exponential with parameter 1.
\item[(b)]  For $j=2,...,d$ the particle jumps from $(\om, x)$ to $(\tau_{{\bf e}_j}\om, x+{\bf e}_j)$ with probability $1/2d$ and to $(\tau_{-{\bf e}_j}\om, x-{\bf e}_j)$ with probability $1/2d $.
\item[(c)] The particle jumps from $(\om, x)$ to $(\tau_{{\bf e}_1}\om, x+{\bf e}_1)$ with probability $1/2d + b(\om)$, and to $(\tau_{-{\bf e}_1}\om, x-{\bf e}_1)$ with probability $1/2d - b(\om)$.
\end{itemize}
If $[\om(t), X(t)] \in \Om \times \Z^d$ is the position of the walk at time $t$ then
\begin{equation} \label{J2}
T_{\eta,\zeta} \vp(\om) = \eta E \left[ \int^\infty_0 dt \ e^{-\eta t} \vp(\om(t)) \exp [-i X(t) \cd \zeta] \ \Big| \  \om(0) = \om, \ X(0) = 0 \right].
\end{equation}
It is clear from the representation (\ref{J2}) that $\| T_{\eta,\zeta}\|_\infty \le 1$.  We conclude from this that (\ref{G2}) is solvable with solution given by
\begin{equation} \label{K2}
\hat v_\ve(\xi,\om) = \hat f_\ve(\xi)\; T_{\ve^2,\ve \xi}(1)(\om), \ \om \in \Om, \ \ \xi \in [-\pi/\ve, \pi/\ve]^d.
\end{equation}
To obtain the homogenization theorem we need then to obtain the limit of the RHS of (\ref{K2}) as $\ve \ra 0$.  To facilitate this we observe from (\ref{H2}) that
\begin{equation} \label{V2}
[\mathcal{L}_\zeta + \eta]1 = 1 + \eta - \frac 1 d \sum^d_{j=1} \cos \zeta_j + 2ib(\om) \sin \zeta_1.
\end{equation}
It follows therefore that
\begin{equation} \label{L2}
T_{\eta,\zeta}(1)(\om) = \eta \Big/ \left[ 1 + \eta - \frac 1 d \sum^d_{j=1} \cos \zeta_j \right]
\end{equation}
\[ 	- \left\{ 2i \sin \zeta_1 \Big/ \left[ 1 + \eta - \frac 1 d \sum^d_{j=1} \cos \zeta_j \right] \right\} T_{\eta,\zeta}b(\om). \]
Setting $\eta = \ve^2, \ \zeta = \ve \; \xi$ for some fixed $\xi \in \R^d$ we see from (\ref{L2}) that 
\begin{equation} \label{M2}
\lim_{\ve \ra 0}T_{\ve^2,\ve\xi}(1)(\om) = 1\Big/ \left[ 1 + \frac {1}{2d} \sum^d_{j=1}  \xi^2_j \right]
\end{equation}
\[ - \left\{ 2i\xi_1 \Big/ \left[ 1 + \frac 1 {2d} \sum^d_{j=1}  \xi^2_j \right] \right\} \lim_{\ve \ra 0}\ve^{-1}
T_{\ve^2,\ve\xi}b(\om).    \]
We shall show that under the assumption of (\ref{C1}) the limit on the RHS of (\ref{M2}) exists.  To do this we define two subspaces of the space $L^\infty(\Om)$.  We define $L^\infty_R(\Om)$ as all functions $\Phi \in L^\infty(\Om)$ such that
\[ 
\left< \Phi(\tau_x \cd) \prod^n_{i=1} b(\tau_{x_i} \cd) \right> = (-1)^{n+1} 
\left<  \Phi(\tau_{Rx}\cdot)  \prod^n_{i=1} b(\tau_{Rx_i} \cd) \right>, \  x, x_i \in \Z^d, \ 1 \le i \le n, n=0,1,2,...
\]
Evidently (\ref{C1}) implies that $b \in L^\infty_R(\Om)$.  We also see that if $\Phi \in L^\infty_R(\Om)$ then $\Phi(\tau_{\bf e_j }\cd)$ and $\Phi(\tau_{-{\bf e_j}} \cd)$ are also in $L^\infty_R(\Om), j=2,...,d$.  For $j=1$ one has that  $\Phi \in L^\infty_R(\Om)$ implies both $[\Phi(\tau_{\bf e_1}  \cd)+ \Phi(\tau_{-{\bf e_1}} \cd)]$ and $b(\cd) [\Phi(\tau_{\bf e_1} \cd)- \Phi(\tau_{-{\bf e_1}} \cd)]$ are in $L^\infty_R(\Om)$.  The space $\hat L^\infty_R(\Om)$ is defined similarly as all functions $\Phi \in L^\infty(\Om)$ such that
\[	
\left< \Phi(\tau_x \cd) \prod^n_{i=1} b(\tau_{x_i} \cd) \right> = (-1)^{n} \left<  \Phi(\tau_{Rx} \cd) \prod^n_{i=1} b(\tau_{Rx_i} \cd) \right>,  x, x_i \in \Z^d, 1 \le i \le n, n = 0,1,2,... .
\]
From (\ref{C1}) we see that the function $\Phi \equiv 1$ is in $\hat L^\infty_R(\Om)$.  As for the space $L^\infty_R(\Om)$, if $\Phi \in \hat L^\infty_R(\Om)$ then 
$\Phi(\tau_{\bf e_j} \cd)$ and $\Phi(\tau_{-{\bf e}_j} \cd)$ are also in $\hat L^\infty_R(\Om), j=2,...,d$.  Similarly both $[\Phi(\tau_{\bf e_1} \cd)+ \Phi(\tau_{-{\bf e}_1} \cd)]$ and $b(\cd) [\Phi(\tau_{\bf e_1} \cd)- \Phi(\tau_{-{\bf e}_1} \cd)]$ are in $\hat L^\infty_R(\Om)$.  We note that the mapping $\Phi(\cd) \ra b(\cd) \Phi(\cd)$ maps $L^\infty_R(\Om)$ into $\hat L^\infty_R(\Om)$ and vice-versa.

We denote the operator $\mathcal{L}_\zeta$ of (\ref{H2}) for $\zeta = 0$ by $\mathcal{L}$.  It is evident that $\mathcal{L}$ is the generator of a random walk on $\Om$.  Hence the kernel of the operator $\mathcal{L}$ is just the constant function.  Furthermore $\mathcal{L}$ leaves the space $L^\infty_R(\Om)$ invariant.  Since the constant function is not in $L^\infty_R(\Om)$ it follows that there is a unique function $\vp \in L^\infty_R(\Om)$ such that
\begin{equation} \label{N2}
\mathcal{L}\vp = b.
\end{equation}
Let $\vp^*$ be the invariant measure for the walk on $\Om$ generated by $\mathcal{L}$.  Thus $\vp^* > 0$,
\begin{equation} \label{O2}
\mathcal{L}^* \vp^* = 0, \ \ \ \ \ \ \left< \vp^* \right> = 1,
\end{equation}
where $\mathcal{L}^*$ is the adjoint of $\mathcal{L}$.  Since $\mathcal{L}$ is non singular on the space $L^\infty_R(\Om)$ it follows that  $\vp^*$ is orthogonal to $L^\infty_R(\Om)$.  We can also see that $\vp^* \in \hat L^\infty_R(\Om)$.  One simply notes that both $\mathcal{L}$ and $\mathcal{L}^*$ leave the space $\hat L^\infty_R(\Om)$ invariant and that the constant function is in $\hat L^\infty_R(\Om)$.  We obtain the limit on the RHS of (\ref{M2}) in terms of the functions $\vp, \vp^*$ defined by (\ref{N2}), (\ref{O2}).
\begin{lemma} Let $\psi \in L^\infty(\Om)$ be defined by
\begin{equation} \label{W2}
\psi(\cd) = \left\{ \frac 1{2d} + b(\cd) \right\} \vp(\tau_{{\bf e}_1} \cd) - \left\{ \frac 1{2d} - b(\cd) \right\} \vp(\tau_{-{\bf e}_1} \cd),
\end{equation}
where $\vp$ is given by (\ref{N2}).  Then if $\vp^*$ is as in (\ref{O2}) there is the limit,
\begin{equation} \label{X2}
\lim_{\ve \ra 0}\ve^{-1} T_{\ve^2,\ve\xi}b(\om) = -i \xi_1 <\vp^*\psi>\Big/ \left[ 1 +\frac 1 {2d} \sum^d_{j=1} \xi^2_j + 2\xi^2_1 <\vp^*\psi> \right],
\end{equation}
for all $\om \in \Om$, provided $\xi_1$ is sufficiently small.
\end{lemma}
\begin{proof}  For $\eta > 0, \ \zeta \in [-\pi, \pi]^d$, let $\vp(\eta,\zeta)$ be the unique solution to the equation
\begin{equation} \label{P2}
[\mathcal{L}_\zeta + \eta] \vp(\eta, \zeta) = b.
\end{equation}
It is clear then that
\begin{equation} \label{Q2}
\ve^{-1} T_{\ve^2,\ve\xi} \ b = \ve \vp(\ve^2, \ve \xi).
\end{equation}
We define operators $A_\zeta, B_\zeta$ by $A_\zeta = [\mathcal{L}_\zeta + \mathcal{L}_{R\zeta}]/2$, $B_\zeta = [\mathcal{L}_\zeta - \mathcal{L}_{R\zeta}]/2$, where $R(\zeta_1, ..., \zeta_d) = (-\zeta_1, \zeta_2,...,\zeta_d)$.  The operator $A_\zeta$ leaves the spaces $L^\infty_R(\Om), \hat L^\infty_R(\Om)$ invariant whereas $B_\zeta$ takes $L^\infty_R(\Om)$ into $\hat L^\infty_R(\Om)$ and vice versa.  Equation (\ref{P2}) is now equivalent to
\begin{equation} \label{T2}
[(A_\zeta + \eta) + B_\zeta] \ \vp(\eta,\zeta) = b \ ,
\end{equation}
and we may write the solution of this formally as a power series,
\begin{equation} \label{R2}
\vp(\eta,\zeta) = \sum^\infty_{n=0} \left\{ -(A_\zeta + \eta)^{-1} B_\zeta\right\}^n \ (A_\zeta + \eta)^{-1} \ b.
\end{equation}
The operators $B_\zeta$, $(A_\zeta + \eta)^{-1} $ on $L^\infty(\Om)$ have norms satisfying $\| B_\zeta\| \le C_2|\zeta_1|$, $\| (A_\zeta + \eta)^{-1} \| \le 1/\eta$, for some constant $C_2$.  Since $A_0 = \mathcal{L}$ is invertible on $L^\infty_R(\Om)$ it follows that for $(\eta, \zeta)$ sufficiently small the operator norm of $(A_\zeta + \eta)^{-1} $ acting on $L^\infty_R(\Om)$ satisfies $\| (A_\zeta + \eta)^{-1} \| \le C_1$ for some constant $C_1$.  We conclude therefore for $(\eta,\zeta)$ sufficiently small that
\begin{equation} \label{S2}
\left\| \left\{ (A_\zeta + \eta)^{-1}B_\zeta \right\}^n (A_\zeta + \eta)^{-1}b\right\|_\infty \le C^n_2 |\zeta_1|^n \; C^{n+1-r}_1 \eta^{-r} \|b\|_\infty \ , 
\end{equation}
where $r = n/2$ if $n$ is even, $r = (n+1)/2$ if $n$ is odd.  Hence if $(\eta, \zeta)$ and $|\zeta_1|^2/\eta$ are small then the series in (\ref{R2}) converges in $L^\infty(\Om)$ to the solution of $(\ref{T2})$.  It follows that for $\xi \in \R^d$ fixed with $\xi_1$ sufficiently small we may construct the function $\vp(\ve^2, \ve \xi)$ by means of (\ref{R2}) as $\ve \ra 0$.

To obtain the limit in (\ref{X2}) we write
\begin{equation} \label{AF2}
\vp(\eta, \zeta) = \vp_1(\eta, \zeta) + \vp_2(\eta, \zeta)
\end{equation}
where $\vp_1(\eta, \zeta) $ is the sum on the RHS of (\ref{R2}) over odd powers of $n$.  It is evident from (\ref{S2}) that for $|\xi_1| < 1 \big/ C_2\sqrt{{C_1}}$ one has
\begin{equation} \label{Y2}
\lim_{\ve \ra 0} \ \ve \vp_2(\ve^2, \ve \xi) = 0.
\end{equation}
We consider the first term in the sum for $\vp_1$.  Setting $\eta = \ve^2, \ \zeta = \ve \xi$ and multiplying the term by $\ve$ as in (\ref{Q2}) we see that
$$
-\lim_{\ve \ra 0}  \ve \left(A_{\ve\xi} + \ve^2\right)^{-1} B_{\ve\xi} \left(A_{\ve\xi} + \ve^2\right)^{-1} b 
$$
$$
=- \lim_{\ve \ra 0} i \ve^2 \xi_1 \left(A_{\ve\xi} + \ve^2\right)^{-1} \bigg[ \bigg\{ \frac 1{2d} + b(\cd)\bigg\} \vp(\tau_{{\bf e}_1}\cd) -  \bigg\{ \frac 1{2d} - b(\cd)\bigg\} \vp(\tau_{-{\bf e}_1}\cd) \bigg],
$$
where $\vp$ is the solution of (\ref{N2}).  Observe now that for any $\psi \in L^\infty(\Om)$ we have 
\[- \lim_{\ve \ra 0} i \ve^2 \xi_1 \left(A_{\ve\xi} + \ve^2\right)^{-1} \psi = -i\xi_1\big<\vp^* \psi\big> 
\lim_{\ve \ra 0}  \ve^2 \left(A_{\ve\xi} + \ve^2\right)^{-1} \; 1 \ ,   \]
where $\vp^*$ is the solution of (\ref{O2}).  From (\ref{V2}) we have that
\[\ve^2 \left(A_{\ve\xi} + \ve^2\right)^{-1} \; 1 = \ve^2 \bigg/ \left[ 1 + \ve^2 - \frac 1 d \sum^d_{j=1} \cos \ve \xi_j \right],\]
whence we conclude that
\[ \lim_{\ve \ra 0}  \ve^2 \left(A_{\ve\xi} + \ve^2\right)^{-1} \; 1 = 1 \bigg/ \left[ 1 + \frac{1}{2d} \sum^d_{j=1}  \xi^2_j \right].
\]
We have therefore obtained a formula for the limit as $\ve \ra 0$ of the first term in the series representation of $\ve \vp_1(\ve^2, \ve \xi)$.  Using the same argument we can obtain a formula for all the terms.  For the rth term corresponding to $r = (n+1)/2$ with $n$ as in (\ref{R2}) we see that the limit is given by the formula,
\begin{equation} \label{AB2}
\frac i{2\xi_1} \left\{ -2 \xi^2_1  \big<\vp^* \psi\big> \bigg/ \left[ 1 + \frac{1}{2d} \sum^d_{j=1}  \xi^2_j \right]\right\}^r,
\end{equation}
where $\psi$ is the function (\ref{W2}).  Evidently $\psi \in \hat L^\infty_R(\Om)$.  We have already observed that $\vp^*$ is also in $ \hat L^\infty_R(\Om)$.  We conclude that
\[
\lim_{\ve \ra 0} \ve \vp_1(\ve^2, \ve \xi) = -i\xi_1 \big<\vp^* \psi\big>  \bigg/ \left[ 1 + \frac{1}{2d} \sum^d_{j=1}  \xi^2_j + 2\xi^2_1 \big<\vp^* \psi\big> \right].
\]
Then (\ref{X2}) follows from this and (\ref{Y2}).  
\end{proof}
Lemma 2.1 enables us to compute the limit (\ref{M2}) when $\xi_1$ is small.  We have
\begin{equation} \label{Z2}
\lim_{\ve \ra 0} T_{\ve^2,\ve\xi}(1)(\om) = 1\bigg/ \left[ 1 + \frac{1}{2d} \sum^d_{j=1}  \xi^2_j + 2\xi^2_1 \big<\vp^* \psi\big> \right].
\end{equation}
We wish now to extend the identity (\ref{Z2}) to all $\xi \in \R^d$.
\begin{lemma} Let $K \subset \R^d$ be a compact set.  Then the limit (\ref{X2}) is uniform for $\xi \in K, \ \om \in \Om$.
\end{lemma}
\begin{proof}  Since the LHS of (\ref{Z2}) does not exceed 1 in absolute value we conclude that
\begin{equation} \label{AA2}
\frac 1{2d} + 2 \big<\vp^* \psi\big>  \ge 0.
\end{equation}
The inequality  (\ref{AA2}) in turn implies that the expression (\ref{AB2}) is the rth power of a number strictly less than 1 provided we also assume that
\begin{equation} \label{AC2}
2 \; \big<\vp^* \psi\big> \le 1/2d.
\end{equation}
We show that the power series methods of Lemma 2.1 apply to prove the result under the additional assumption (\ref{AC2}).  We shall see in $\S$3 that  $\big<\vp^* \psi\big> \le 0$ for dimension $d=1$, in which case (\ref{AC2}) certainly holds.    Since the constant function is the unique eigenvector of $A_0 = \mathcal{L}$ with eigenvalue 0 and it is also an eigenvector of $A_\zeta$ it follows that there exists $\de > 0$ such that if $|\zeta| < \de$ then the adjoint $A^*_\zeta$ of $A_\zeta$  has a unique eigenvector $\vp^*_\zeta$ with eigenvalue equal to the eigenvalue of $A_\zeta$ for the constant function.  Normalizing $\vp^*_\zeta$ so that $<\vp^*_\zeta> = 1$, it is easy to see that there is a constant $C_1$ such that
\begin{equation} \label{AD2}
\| \vp^*_\zeta - \vp^*\| _\infty \le C_1|\zeta|, \ \ |\zeta| < \de.
\end{equation}
For $|\zeta| < \de, \ \eta > 0$ we define a projection $P_{\eta,\zeta}$ by 
\[	P_{\eta,\zeta} \; \psi = \left< \bar \vp^*_\zeta \; \psi \right> \left( A_\zeta + \eta\right)^{-1} \ 1, \ \psi \in L^\infty(\Om). \]
Then there is a constant $C_2$ such that 
\begin{equation} \label{AE2}
\| (A_\zeta + \eta)^{-1} -P_{\eta,\zeta}\| \le C_2, \ |\zeta| < \de.
\end{equation}
The uniform convergence of (\ref{X2}) for $\xi \in K$ follows now from(\ref{AD2}), (\ref{AE2}) just as in Lemma 2.1.

Finally we consider the situation where (\ref{AC2}) is violated.  As in Lemma 2.1 we decompose the solution $\vp(\eta,\zeta)$ of (\ref{P2}) into a sum (\ref{AF2}).  The function $\vp_1(\eta, \zeta)$ is a solution to the equation,
\begin{equation} \label{AG2}
\left[ A_\zeta + \eta - B_\zeta \left( A_\zeta + \eta\right)^{-1} B_\zeta \right] \vp_1(\eta, \zeta) = -B_\zeta (A_\zeta + \eta)^{-1} \; b.
\end{equation}
The function $\vp_2(\eta, \zeta)$ is a solution to the equation,
\begin{equation} \label{AH2}
\left[ A_\zeta + \eta - B_\zeta \left( A_\zeta + \eta\right)^{-1} B_\zeta \right] \vp_2(\eta, \zeta) = b.
\end{equation}
It is easy to see that if $\vp_2(\eta, \zeta)$ is a solution of  (\ref{AH2}) then the function $\vp_1(\eta, \zeta) = \vp(\eta, \zeta) - \vp_2(\eta, \zeta)$, where $\vp(\eta, \zeta)$ solves (\ref{P2}), is a solution to (\ref{AG2}).  Hence if (\ref{AG2}),(\ref{AH2}) have unique solutions $\vp_1(\eta, \zeta), \vp_2(\eta, \zeta)$ then the identity (\ref{AF2}) holds.

We show that (\ref{AH2}) has a unique solution in $L^\infty_R(\Om)$ provided $\eta > 0$ and $(\eta, \zeta)$ are sufficiently small.  To see this we write (\ref{AH2}) as
\begin{equation} \label{AI2}
\left[ A_\zeta + \eta - L_{\eta,\zeta} - B_\zeta  P_{\eta,\zeta} B_\zeta \right] \vp_2(\eta, \zeta) = b,
\end{equation}
where by (\ref{AE2}) the operator $L_{\eta,\zeta}$ is invariant on $L^\infty_R(\Om)$ and satisfies $\| L_{\eta,\zeta}\| \le C|\zeta|^2$ for some constant $C$.  Next let $\vp_3(\eta,\zeta)$ be the solution to
\begin{equation} \label{AJ2}
\left[ A_\zeta + \eta - L_{\eta,\zeta}  \right] \vp_3(\eta, \zeta) = b.
\end{equation}
For $(\eta,\zeta)$ small there is a unique solution to (\ref{AJ2}) in $L^\infty_R(\Om)$ which satisfies
\begin{equation} \label{AK2}
 \|\vp - \vp_3(\eta,\zeta)\|_\infty \le C[|\eta| + |\zeta|^2],
 \end{equation}
 for some constant $C$, where $\vp$ is the solution of (\ref{N2}).  Now it is easy to see that the solution $\vp_2(\eta, \zeta)$ of (\ref{AI2}) is given in terms of $\vp_3(\eta,\zeta)$ by the formula,
\begin{multline} \label{AL2}
\vp_2(\eta, \zeta) = \left[ 1 + \eta - \frac 1 d \sum^d_{j=1} \cos \zeta_j \right] \vp_3(\eta, \zeta) \Big/
\\
\left[ 1 + \eta - \frac 1 d \sum^d_{j=1} \cos \zeta_j  + 2 \left< \bar{\vp}^*_\zeta B\vp_3(\eta, \zeta) \right> \sin^2 \zeta_1\right], 
\end{multline}
where the operator $B$ is defined by $B_\zeta = i \sin \zeta_1 B$.  In view of (\ref{AK2}) and the fact that $B\vp = \psi$ and we are assuming (\ref{AC2}) is violated, it follows that the denominator in (\ref{AL2}) is positive for $(\eta, \zeta)$ sufficiently small.  We have shown a solution $\vp_2(\eta, \zeta)$ of (\ref{AI2}) exists in $L^\infty_R(\Om)$. The uniqueness of the solution follows from the uniqueness of the solution to (\ref{AJ2}).  Evidently the limit (\ref{Y2}) follows from (\ref{AK2}), (\ref{AL2}) for all $\xi$ and is uniform for $\xi$ restricted to a compact subset of $\R^d$.

Next we show that (\ref{AG2}) has a unique solution in $\hat L^\infty_R(\Om)$ provided $\eta > 0$ and $(\eta,\zeta)$ are sufficiently small.  First note that for $(\eta, \zeta)$ small the operator  $B_\zeta(A_\zeta + \eta)^{-1} B_\zeta$ leaves $\hat L^\infty_R(\Om)$ invariant and there is a constant $C$ such that
\begin{equation} \label{AM2}
\| B_\zeta(A_\zeta + \eta)^{-1} B_\zeta\| \le C|\zeta|^2.  
\end{equation}

We define the subspace $\mathcal{E}_\zeta$ of $\hat L^\infty_R(\Om)$ by
\[	\mathcal{E}_\zeta = \left\{ \psi \in \hat L^\infty_R(\Om): \left< \psi \  \bar \vp^*_\zeta \right> = 0 \right\} . \]
Let $P_\zeta$ be the projection operator on $\hat L^\infty_R(\Om)$ orthogonal to $\mathcal{E}_\zeta$, whence
\[	P_\zeta \psi = \left< \psi  \ \bar \vp^*_\zeta \right>,  \quad \psi \in \hat L^\infty_R(\Om). \]
Consider now the equation related to (\ref{AG2}) given by
\begin{eqnarray} \label{AN2}
\bigg[  A_\zeta + \eta &-& (I-P_\zeta)B_\ze (A_\zeta + \eta)^{-1} B_\ze \bigg] \vp_4(\eta,\ze) \\
&=& -(I-P_\zeta)B_\ze (A_\zeta + \eta)^{-1} b. \nn
\end{eqnarray}
In view of (\ref{AM2}) it is clear that for $(\eta,\ze)$ sufficiently small the equation (\ref{AN2}) has a unique solution $\vp_4(\eta,\ze)$ in $\mathcal{E}_\ze$.  Furthermore, if we define $\vp_1(\eta,\ze)$ by
\begin{multline} \label{AO2}
\vp_1(\eta,\ze) = \bigg\{ \left[ 1 + \eta - \frac 1 d \sum^d_{j=1} \cos \zeta_j \right] \vp_4(\eta,\ze) 
- i \sin \ze_1 \left< \bar \vp^*_\ze B(A_\ze + \eta)^{-1} \; b \right> 
\\
- \sin^2 \ze_1 \left< \bar \vp^*_\ze B(A_\ze + \eta)^{-1} \; B \; \vp_4(\eta,\ze)\right> \bigg\} \Bigg{/}  
\\
 \bigg\{ 1 + \eta - \frac 1 d \sum^d_{j=1} \cos \zeta_j + 2\sin^2 \ze_1 \left< \bar\vp^*_\ze B(A_\ze + \eta)^{-1} \; b \right> 
 \\
 -2i \sin^3 \ze_1 \left< \bar\vp^*_\ze B(A_\ze + \eta)^{-1} \; B  \vp_4 \right> \bigg\} , 
 \end{multline}
then one sees that the formula (\ref{AO2}) yields a solution to (\ref{AG2}).  Conversely, since we are assuming (\ref{AC2}) is violated, it follows that for $(\eta, \zeta)$ small (\ref{AO2}) is the unique solution in $\hat L^\infty_R(\Om)$ to (\ref{AG2}).  It is easy to see now from (\ref{AO2}) that the limit $\lim_{\ve \ra 0} \ve\vp_1(\ve^2,\ve\xi)$ exists and is uniform for $\xi$ in a compact subset of $\R^d$.  Furthermore, the limit is given by the RHS of (\ref{X2}).

Finally we show that if $\vp_1(\eta,\zeta), \vp_2(\eta,\zeta)$ are solutions to (\ref{AG2}), (\ref{AH2}) then (\ref{AF2}) holds.  To see this we put $\vp(\eta,\zeta) =\vp_1(\eta,\zeta) $ $+ \vp_2(\eta,\zeta)$ and note that (\ref{AG2}), (\ref{AH2}) imply that $\vp(\eta,\zeta)$ satisfies the equation
\[  \left[  A_\zeta + \eta - B_\zeta (A_\zeta + \eta)^{-1} B_\ze \right] \vp(\eta,\ze) 
= b - B_\ze (A_\zeta + \eta)^{-1} b. \]
We can rewrite this equation as
\[ \left[  A_\zeta +  B_\ze + \eta - B_\zeta (A_\zeta + \eta)^{-1}(A_\ze+ B_\ze+\eta) \right] \vp(\eta,\ze) 
= b - B_\ze (A_\zeta + \eta)^{-1} b, \]
which is the same as 
\[ 	[\mathcal{L}_{R\zeta} + \eta ](A_\ze + \eta)^{-1} [\mathcal{L}_\ze + \eta] \vp(\eta, \ze) = [\mathcal{L}_{R\zeta} + \eta ](A_\ze + \eta)^{-1} \; b.  \]
Now using the fact that the operator $\mathcal{L}_{R\zeta} + \eta$ is invertible we obtain (\ref{P2}). 
\end{proof}
Next we show that there is strict inequality in (\ref{AA2}).  In order to do this we shall first obtain a concrete representation of spaces $\Om$ which satisfy (\ref{B1}).
\begin{lemma}  Let $\Om$ be a finite probability space and $b : \Om \ra \R$ satisfy (\ref{C1}).  Then $\Om$ may be identified with a rectangle in $\Z^d$ with periodic boundary conditions.  The operators $\tau_x, \; x \in \Z^d$, act on  $\Om$ by translation and the measure $\left< \cd \right>$ is simple averaging.  Let $R : \Om \ra \Om$ be the reflection operator defined as reflection in the hyperplane through the center of $\Om$ with normal ${\bf e}_1$.  Then there is the identity $b(\om) = -b(R\om), \ \om \in \Om$.
\end{lemma}
\begin{proof}
Since $\Om$ has no nontrivial invariant subsets under the action of the  $\tau_{{\bf e}_j}, \; 1 \le j \le d$, it is isomorphic to a rectangle in $\Z^d$ with periodic boundary conditions.  Thus we may assume $\Om$ is given by
\begin{equation} \label{AP2}
\Om = \left\{ x = (x_1, ..., x_d) \in \Z^d : 0 \le x_i \le L_i - 1, \ 1 \le i \le d \right\},
\end{equation}
where $L_1,...,L_d$ are positive integers.  The action of the $\tau_{{\bf e}_j}$ is translation, $\tau_{{\bf e}_j} x = x + {\bf e}_j$ with periodic boundary conditions.  The measure on $\Om$ is averaging,
\begin{equation} \label{AQ2}
\left< \Psi(\cd) \right> = \frac 1{L_1L_2 \cdots L_d} \ 
\begin{array}[t]{c}
{\displaystyle\sum_{0\le x_i \le L_i-1,}}\\
{\scriptstyle 1 \le i\le d}
\end{array} 
\Psi(x_1,...,x_d).
\end{equation}
Functions $\Psi: \Om \ra \R$ are isomorphic to periodic functions $\Psi : \Z^d \ra \R$.

Next we consider the condition  (\ref{C1}).  We define a function $b^R : \Om \ra \R$ by $b^R(\om) = -b(R\om)$,  $\om \in \Om$.  It is easy to see that $b^R(\tau_x \om) = -b(R\tau_x\om)$ $= -b(\tau_{Rx}  R\om), \om \in \Om$.  Since $R$ leaves the measure (\ref{AQ2}) invariant (\ref{C1}) implies that for any $\theta_1,...,\theta_k \in \R$, $x_1,...,x_k \in \Z^d$, there is the identity,
\[
\left< \exp \left[ \sum^k_{j=1} \theta_j \ b(\tau_{x_j} \cd ) \right] \right> = 
\left< \exp \left[ \sum^k_{j=1} \theta_j \ b^R(\tau_{x_j} \cd) \right] \right>.
\]
We conclude that $b \equiv b^R$. 
\end{proof}
Next we wish to construct the solutions $\vp, \vp^*$ of (\ref{N2}), (\ref{O2}) on the domain $\Om$ defined by (\ref{AP2}).  First observe that since $\Om$ is the fundamental region for the homogenization problem we can assume that $L_1$ is an even integer by simply doubling $\Om$ if $\L_1$ is odd.  In that case the function $b$ is determined by its values $b(x), \ x \in \Om$, $0 \le x_1 \le L_1/2 - 1$.  Hence we define a new fundamental region $\hat \Om$ by
\begin{equation} \label{AR2}
\hat \Om = \{ x \in \Om : 0 \le x_1 \le L_1/2 - 1 \}.
\end{equation}
We can extend functions $\Psi : \hat \Om \ra \R$ to $\Om$  by either symmetric or antisymmetric extension.  For a symmetric extension we define $\Psi$ on $\Om - \hat \Om$ by
\begin{equation} \label{AS2}
\Psi(x_1,...,x_d) = \Psi(L_1 - 1 - x_1, x_2, ..., x_d), \ L_1/2 \le x_1 \le L_1 - 1.
\end{equation}
For an antisymmetric extension we define $\Psi$ by
\begin{equation} \label{AT2}
\Psi(x_1,...,x_d) = -\Psi(L_1 - 1 - x_1, x_2, ..., x_d), \ L_1/2 \le x_1 \le L_1 - 1.
\end{equation}
\begin{lemma} The solution $\vp : \Om \ra \R$ of (\ref{N2}) is an antisymmetric extension of its restriction to $\hat \Om$.  The solution $\vp^* : \Om \ra \R$  of (\ref{O2}) is a symmetric extension of its restriction to $\hat \Om$.
\end{lemma}
\begin{proof}
This follows easily from the fact that $b : \Om \ra \R$ is an antisymmetric extension of its restriction to $\hat \Om$ and the uniqueness of the solution to (\ref{N2}), (\ref{O2}). 
\end{proof}
Lemma 2.4 implies that we can find the functions $\vp, \vp^*$ by solving (\ref{N2}), (\ref{O2}) on $\hat \Om$ with antisymmetric and symmetric boundary conditions respectively.  Thus $\mathcal{L}$ acting on functions $\Psi : \hat \Om \ra \R$ with antisymmetric boundary conditions is defined by 
\begin{multline} \label{AU2}
\mathcal{L} \Psi(x) = \Psi(x) - \sum^d_{j=1} \frac 1{2d} \left[ \Psi(x + {\bf e}_j) + \Psi(x - {\bf e}_j) \right] \\
- b(x) \left[ \Psi(x + {\bf e}_1) - \Psi(x - {\bf e}_1) \right] , \ x \in \hat \Om, 
\end{multline}
where the boundary conditions are given by,
\begin{multline} \label{AV2}
\Psi(-1, x_2,...,x_d) = -\Psi(0, x_2,...,x_d),  \quad
\Psi(L_1/2, x_2,...,x_d) = -\Psi(L_1/2 - 1, x_2,...,x_d) ,  \\
0 \le x_j \le L_j - 1, \ j=2,...,d, 
\end{multline}
and periodic boundary conditions in the directions ${\bf e}_j, \ 2 \le j \le d$.  Evidently (\ref{AV2}) is derived from (\ref{AT2}).  It is easy to see that the operator $\mathcal{L}$ is invertible on the space $L^\infty(\hat \Om)$ if the boundary conditions (\ref{AV2}) are imposed.  In fact the solution to the equation
\begin{equation} \label{AW2}
\mathcal{L} \; \Psi(x) = f(x), \ \ x \in \hat \Om,
\end{equation}
with boundary conditions (\ref{AV2}) can be represented as an expectation for a continuous time Markov chain $X(t), t \ge 0$, on $\hat \Om$.  For the chain the transition probabilities at a site $x \in \hat \Om$ satisfying $0 < x_1 < L_1/2 - 1$ are given by $x \ra x + {\bf e}_j, x \ra x - {\bf e}_j$,  $2 \le j \le d$, each with probability $1/2d$, with periodic boundary conditions in direction ${\bf e}_j, 2 \le j \le d$.  In the direction ${\bf e}_1$ then $x \ra x + {\bf e}_1$ with probability $1/2d + b(x)$ and $x \ra x - {\bf e}_1$ with probability $1/2d - b(x)$.  The waiting time at site $x$ is exponential with parameter 1.  If $x_1 = 0$ then $x \ra x \pm {\bf e}_j$, $2 \le j \le d$, with probability $1/2d[1+1/2d - b(x)] < 1/2d$, and $x\ra x + {\bf e}_1$ with probability $[1/2d + b(x)] / [1+1/2d - b(x)] < 1/d$.  The waiting time is exponential with parameter $ [1+1/2d - b(x)]$.  Note that there is a positive probability that the walk will be killed at a site $x$ with $x_1 = 0$.  A similar situation occurs at a site $x$ with $x_1 = L_1/2 - 1$.  Now $x \ra x - {\bf e}_1$ with probability $[1/2d  - b(x)] / [1+1/2d + b(x)] < 1/d$ and the waiting time is exponential with parameter $[1+1/2d + b(x)].$ The solution $\Psi$ of (\ref{AW2}) with boundary conditions (\ref{AV2}) has the representation
\begin{equation} \label{BA2}
\Psi(x) = E \left[ \int^\tau_0  f (X(t))dt | X(0) = x \right], \ \ x \in \hat \Om,
\end{equation}
where $\tau$ is the killing time for the chain. 

We may also consider the operator $\mathcal{L}$ of (\ref{AU2}) with symmetric boundary conditions,
\begin{multline} \label{AX2}
\Psi(-1, x_2,...,x_d) = \Psi(0, x_2,...,x_d), \quad
\Psi(L_1/2, x_2,...,x_d) = \Psi(L_1/2 - 1, x_2,...,x_d) , \\ 0 \le x_j \le L_j - 1, \ j=2,...,d, 
\end{multline}
corresponding to (\ref{AS2}).  This is also associated with a continuous time Markov chain $X(t)$ on $\hat \Om$.  The transition probabilities and waiting time at a site $x \in \hat \Om$ with $0 < x_1 < L_1/2 - 1$ are as for the chain defined in the previous paragraph.  For $x \in \hat \Om$ with $x_1=0$ reflecting boundary conditions corresponding to (\ref{AX2}) are imposed.  Thus the waiting time at $x$ is exponential with parameter $[1 - 1/2d + b(x)], x \ra x \pm {\bf e}_j$, $2 \le j \le d$, with probability $1/2d[1-1/2d + b(x)]$ and $x\ra x+{\bf e}_1$ with probability $[1/2d+b(x)] / [1-1/2d + b(x)]$.  A similar situation occurs at $x \in \hat \Om$ with $x_1 = L_1/2 - 1$.  The formal adjoint $\mathcal{L}^*$ of the operator $\mathcal{L}$ of (\ref{AU2}) is given by
\begin{eqnarray} \label{AY2}
\mathcal{L}^* \Psi(x) &=& \Psi(x) - \sum^d_{j=1} \frac 1{2d} \left[ \Psi(x + {\bf e}_j) + \Psi(x - {\bf e}_j) \right] \\
&-& b(x-{\bf e}_1)  \Psi(x - {\bf e}_1) + b(x + {\bf e}_1) \Psi(x + {\bf e}_1), \ x \in \hat \Om. \nn
\end{eqnarray}
It is easy to see that for functions $\Phi, \Psi$ on $\hat \Om$ satisfying the symmetric boundary conditions (\ref{AX2}) there is the identity
\begin{equation} \label{AZ2}
\left< \Phi \; \mathcal{L}^* \; \Psi \right>_{\hat \Om} = \left< \Psi \; \mathcal{L} \; \Phi\; \right>_{\hat \Om},
\end{equation}
where $\left< \cd \right>_{\hat \Om}$ is the uniform probability measure on $\hat \Om$.  Note that to show (\ref{AZ2})  one has to use the fact that the function $b$ satisfies the antisymmetric conditions (\ref{AV2}).  Hence the adjoint of the operator $\mathcal{L}$ acting on functions $\Psi : \hat \Om \ra \R$ with symmetric boundary conditions (\ref{AX2}) is the operator $\mathcal{L}^*$ of (\ref{AY2}) also acting on functions with symmetric boundary conditions.  In particular, it follows from Lemma 2.4 that the solution $\vp^*$ of (\ref{O2}), restricted to $\hat \Om$, is the unique invariant measure for the Markov chain $X(t)$.

Next let $\psi_0 : \hat \Om \ra \R$ be the solution of the homogeneous equation (\ref{AW2}) i.e. $f \equiv 0$, with the non-homogeneous antisymmetric boundary conditions
\begin{multline} \label{BB2}
\psi_0(-1, x_2,...,x_d) = -\psi_0(0, x_2,...,x_d), \quad
\psi_0(L_1/2, x_2,...,x_d) = 1 -\psi_0(L_1/2 - 1, x_2,...,x_d) , \\ 0 \le x_j \le L_j - 1, \ j=2,...,d. 
\end{multline}
One can see that $\psi_0$ is a positive function since it has a representation given by (\ref{BA2}), where $f$ is the function
\begin{eqnarray*}
f(x) &=& \frac 1{2d} + b(x), \ x \in \hat \Om, \ \ x_1 = L_1/2 - 1, \\
&=& \ 0, \ \ \ {\rm otherwise}.
\end{eqnarray*}
The following lemma now shows that there is strict inequality in (\ref{AA2})
\begin{lemma} Let $\vp^*$ be the solution of (\ref{O2}) and $\psi$ be given by (\ref{W2}).  Then there is the identity,
\begin{equation} \label{BC2}
\frac 1{2d} + 2 \left<\vp^* \psi\right> = L^2_1 \left< \vp^*(\cd) \left[ \frac 1{2d} - b(\cd)\right]\psi_0(\cd) \chi_0(\cd)\right>_{\hat \Om},
\end{equation}
where $\chi_0 : \hat \Om \ra \R$ is defined by $\chi_0(x) = 1$ if $x_1=0, \  \chi_0(x) = 0$, otherwise.
\end{lemma}
\begin{proof}
Since both $\vp^*$ and $\psi$ are symmetric on $\Om$ in the sense of (\ref{AS2}) we may regard them as functions on $\hat \Om$ with symmetric boundary conditions (\ref{AX2}).  We define a function $\psi_1 : \hat \Om \ra \R$ by $\psi_1(x) = [L_1/2 - 1/2 - x_1]\vp(x), \ x \in \hat{\Om}$, where $\vp$ is the solution to (\ref{N2}).  It is easy to see that
\begin{equation} \label{BD2}
\mathcal{L} \psi_1(x) = [L_1/2 - 1/2 - x_1]b(x) + \psi(x), \ \ \ \ x \in \hat \Om, \ 0 < x_1 <L_1/2 - 1.
\end{equation}
We impose now symmetric boundary conditions on $\psi_1$ at $x_1 = 0, x_1 = L_1/2 - 1$.  One sees that (\ref{BD2}) continues to hold at $x_1 = L_1/2 - 1$ but at $x_1 = 0$ there is the formula,
\begin{equation} \label{BE2}
\mathcal{L} \psi_1(x) = [L_1/2 - 1/2]b(x) + \psi(x) - L_1[1/2d - b(x)]\vp(x).
\end{equation}
In deriving (\ref{BE2}) we have used the fact that $\vp$ satisfies antisymmetric boundary conditions at $x_1=0$.  Now from (\ref{O2}), (\ref{AZ2}), (\ref{BD2}), (\ref{BE2}) we have that
\begin{multline}  \label{BF2}
2 \left<\vp^* \psi\right> = 2 \left<\vp^* \psi\right> _{\hat \Om}  =  \\
2L_1 \left< \vp^*(\cd) \left[ \frac 1{2d} - b(\cd)\right] \vp(\cd) \chi_0(\cd) \right>_{\hat \Om} 
- \left< \vp^*(\cd) \left[ L_1 - 1 -2x_1\right] b(\cd) \right>_{\hat \Om}.
\end{multline}

Next we define a function $\psi_2 : \hat \Om \ra \R$ by $\psi_2(x) = (L_1 - 1 -2x_1)^2, x \in \hat \Om$.  Then we have
\begin{equation} \label{BG2}
\mathcal{L}\psi_2(x) = 8[L_1 - 1 -2x_1] b(x) - 4/d, \ x \in \hat \Om, \ 0 < x_1 <L_1/2 - 1.
\end{equation}
Again we impose symmetric boundary conditions on $\psi_2$ at $x_1=0, x_1 = L_1/2 - 1$, in which case (\ref{BG2}) continues to hold at $x_1 = L_1/2 - 1.$ At $x_1=0$ there is the formula
\begin{equation} \label{BH2}
\mathcal{L}\psi_2(x) = 8[L_1 - 1 -2x_1] b(x) - 4/d + 4L_1 [1/2d - b(x)].
\end{equation}
It follows now from (\ref{O2}), (\ref{AZ2}), (\ref{BG2}), (\ref{BH2}) that
\begin{equation} \label{BI2}
-  \left< \vp^*(\cd) \left[ L_1 - 1 -2x_1\right] b(\cd) \right>_{\hat \Om} = - 1/2d + \frac{L_1}2 
\left< \vp^*(\cd) \left[ \frac 1{2d} - b(\cd)\right] \chi_0(\cd) \right>_{\hat \Om},
\end{equation}
where we have used the fact that $\left<\vp^*\right>_{\hat \Om} = 1$.  It follows now from (\ref{BF2}), (\ref{BH2}) that
\begin{equation} \label{BJ2}
1/2d + 2 \left< \vp^*\psi \right> = \frac{L_1}2 \left< \vp^*(\cd) \left[ \frac 1{2d} - b(\cd)\right] [1 + 4\vp(\cd)] \chi_0(\cd) \right>_{\hat \Om} .
\end{equation}

We put now $\psi_0(x) = [2x_1 + 1 + 4\vp(x)]/2L_1$, and it is easy to verify that $\psi_0$ satisfies the homogenous equation (\ref{AW2}) with the boundary conditions (\ref{BB2}).  The result follows then from (\ref{BJ2}). 
\end{proof}
\begin{proof}[Proof of Theorem 1.1]  The proof proceeds identically to the proof of Theorem 1.1 of \cite{cp}, on using lemmas 2.1-2.4. 
\end{proof}
Finally we wish to show that Theorem 1.2 holds to leading order in perturbation theory.
\begin{theorem}  There exists $\de > 0$ such that if $b : \Om \ra \R$ satisfies $0 < \sup_{\om\in\Om} |b(\om)| < \de$ then $q(b) < 1/2d$, provided $d=1$, or $ d > 1$ and $L_1 \le 4$. If $d\ge 2$ and $L_1 \ge 6$ then there exists arbitrarily small $b$ with  $q(b) >1/2d$.
\end{theorem}
\begin{proof}
We shall use the LHS of (\ref{BC2}) as an expression for $q(b)$.  If $b \equiv 0$ then $\vp^* \equiv 1, \vp \equiv 0 \Rightarrow \psi \equiv 0$.  Thus to obtain an expression for $q(b)$ which is correct to second order in perturbation theory we need to expand $\vp^*$ to first order in $b$ and $\vp$ to second order.  We consider first $\vp^*$ which is the solution to (\ref{O2}).  Letting $\De$ be the finite difference Laplacian acting on functions $\Psi : \Om \ra \R$ with periodic boundary conditions,
\[ \De \Psi(x) = \sum^d_{j=1}  \left[ \Psi(x + {\bf e}_j) + \Psi(x - {\bf e}_j) - 2\Psi(x) \right] , \ x \in \Om,
\]
we have from (\ref{AY2}) that (\ref{O2}) is given by
\begin{equation} \label{BK2}
-\frac{\De}{2d} \vp^*(x) + b(x + {\bf e}_1) \vp^*(x + {\bf e}_1) - b(x - {\bf e}_1) \vp^*(x - {\bf e}_1) = 0, \ \  x \in \Om, \  \left<\vp^*\right> = 1.
\end{equation}
Since $\left< \left[ \tau_{-{\bf e}_1} - \tau_{{\bf e}_1} \right] b\right> = 0$ the solution to (\ref{BK2}) is to first order in perturbation theory given by
\begin{equation} \label{BL2}
\vp^* = 1 + (-\De/2d)^{-1} \left[ \tau_{-{\bf e}_1} - \tau_{{\bf e}_1} \right] b.
\end{equation}
From (\ref{AU2}) equation (\ref{N2}) is the same as
\begin{equation} \label{BM2}
-\frac{\De}{2d} \vp(x) - b(x) \left[ \tau_{{\bf e}_1} - \tau_{-{\bf e}_1} \right] \vp(x) = b(x), \ \ x \in \Om.
\end{equation}
Using the fact that
\[
\left<b\right> = \left< b \left[ \tau_{{\bf e}_1} - \tau_{-{\bf e}_1} \right] (-\De/2d)^{-1} b\right> = 0,
\]
we see that the solution to (\ref{BM2}) correct to second order in $b$ is given by
\begin{equation} \label{BN2}
\vp = (-\De/2d)^{-1} b + (-\De/2d)^{-1} b   \left[ \tau_{{\bf e}_1} - \tau_{-{\bf e}_1} \right] (-\De/2d)^{-1} b.
\end{equation}
From (\ref{W2}) and (\ref{BN2}) we can obtain an expression for $\psi$ which is correct to second order in $b$,
\begin{multline} \label{BO2}
\psi = \frac 1{2d} \left[ \tau_{{\bf e}_1} - \tau_{-{\bf e}_1} \right] (-\De/2d)^{-1} b 
+ b \left[ \tau_{{\bf e}_1} + \tau_{-{\bf e}_1} \right] (-\De/2d)^{-1} b  \\
+ \frac 1{2d} \left[ \tau_{{\bf e}_1} - \tau_{-{\bf e}_1} \right] (-\De/2d)^{-1} b \left[ \tau_{{\bf e}_1} - \tau_{-{\bf e}_1} \right] (-\De/2d)^{-1} b. 
\end{multline}
From (\ref{BL2}), (\ref{BO2}) we see that the lowest order term in the expansion of $\left<\vp^* \psi\right>$ in powers of $b$ is second order.  Thus correct to second order we have 
\begin{multline} \label{BP2}
\left<\vp^*\psi\right> = \left< b \left[ \tau_{{\bf e}_1} + \tau_{-{\bf e}_1} \right] (-\De/2d)^{-1} b\right> 
\\
+  \frac 1{2d} \left< b\left[ \tau_{{\bf e}_1} - \tau_{-{\bf e}_1} \right] (-\De/2d)^{-1}  \left[ \tau_{{\bf e}_1} - \tau_{-{\bf e}_1} \right] (-\De/2d)^{-1} b \right>.  
\end{multline}
The RHS of (\ref{BP2}) is a translation invariant quadratic form, whence it has eigenvectors $\exp[i\xi \cd x], \ x \in \Om$, with corresponding eigenvalue given by the formula,
\begin{equation} \label{BQ2}
2d \left\{ \cos \xi_1 - \frac{\sin^2 \xi_1}{\sum^d_{j=1}(1-\cos \xi_j)} \right\} \bigg/ \sum^d_{j=1} (1 - \cos \xi_j).
\end{equation}
We obtain an expression for the quadratic form (\ref{BP2}) by doing an eigenvector decomposition in the $x_1$ direction.  Putting $L = L_1/2$ we have that $\xi_1 = \pi k/L$, $k=0, \pm 1,...,\pm (L-1), L$.  The function $b$ then has a representation,
\[  b(x_1) = \sum_{\xi_1} e^{i\xi_1x_1} \left(  \frac 1{2L} \ \sum^{2L-1}_{y=0} \ b(y)e^{-i\xi_1y} \right).  \]
If we use the antisymmetry property of $b$, $b(y) = -b(2L - 1 - y)$ then one has that
\[   \sum^{2L-1}_{y=0} \ b(y)e^{-i\xi_1y} = -2i e^{i\xi_1/2} \ \sum^{L-1}_{y=0} b(y) \sin \xi_1(y + 1/2).   \]
We conclude from this and (\ref{BQ2}) that the expression (\ref{BP2}) is the same as
\begin{equation} \label{BR2}
\left<\vp^*\psi\right> = -\frac{4d}{L^2} \left< \left( \sum^{L-1}_{y=0} (-1)^y b(y) \right) [-\De_{d-1} + 4]^{-1} \left( \sum^{L-1}_{y=0} (-1)^y b(y) \right)\right>
\end{equation}
\[ + \frac{8d}{L^2} \sum^{L-1}_{k=1} \bigg< \left( \sum^{L-1}_{y=0} b(y) \sin \pi k(y + \frac 1 2)/L \right) \]
\[ \left\{ \cos \frac{\pi k}L - 2 \sin^2 \frac{\pi k}L \left[ -\De_{d-1} + 2 \big( 1 - \cos(\pi k/L) \big)\right]^{-1}\right\}  \]
\[ \left[ -\De_{d-1} + 2 \big( 1 - \cos(\pi k/L) \big)\right]^{-1} \left( \sum^{L-1}_{y=0} b(y) \sin \pi k(y + \frac 1 2)/L \right) \bigg>,	\]
where $\De_{d-1}$ denotes the $d-1$ dimensional Laplacian acting on the space $\{x_1 = 0\}$.

Observe now that the $L$ dimensional vectors $\sin \pi k(y+1/2)/L, 0 \le y \le L-1$, are mutually orthogonal, $k=1,...,L$.  This is a consequence of the fact that they are the eigenvectors of the second difference operator on the set $\{ 0 \le y \le L-1\}$ with antisymmetric boundary conditions.  It follows that the quadratic form (\ref{BR2}) is negative definite if and only if all the eigenvalues (\ref{BQ2}) are negative.  This is the case for $d=1$.  For $d > 1$ it is still true provided $L \le 2$, but already for $L=3$ it is false.  Thus for $L=3$ one can find a $b$ such that the homogenized limit has an effective diffusion constant which is larger than the $b \equiv 0$ case. 
\end{proof}

\section{Proof of Theorem 1.2}
We shall use the representation for the effective diffusion constant given by the RHS of (\ref{BC2}).  We consider first the $d=1$ case.
\begin{lemma} Let $\hat\Om$ be the space $\hat\Om = \{ x \in \Z : 1 \le x \le L\}$.  If $\varphi^* : \Om \ra \R$ is the solution to (\ref{O2}) then $\varphi^*(1)$ is given by the formula,
\begin{equation} \label{A3}
\varphi^*(1)\de_1 = L \; \prod^L_{k=1} \de_k \big/ \sum^L_{r=1} \ \prod^{r-1}_{j=1} \bar\de_j \prod^L_{j=r+1} \de_j,
\end{equation}
where the $\de_j, \bar\de_j, \ 1 \le j \le L$, are given by
\begin{equation} \label{B3}
\de_j = 1/2 - b(j), \ \bar\de_j = 1/2 + b(j).
\end{equation}
\end{lemma}
\begin{proof}
 From (\ref{AY2}) we see that $\varphi^* : \hat\Om \ra \R$ satisfies the equation
\begin{multline}
 \label{C3}
\varphi^*(x) - \frac 1 2 \Big[ \varphi^*(x+1) + \varphi^*(x-1)\Big] - b(x-1) \varphi^*(x-1)\\
+b(x+1)\varphi^*(x+1)=0, \quad 1 \le x \le L,
\end{multline}
with the symmetric boundary conditions and normalization given by
\begin{equation} \label{D3}
\varphi^*(0) = \varphi^*(1), \ \varphi^*(L+1)= \varphi^*(L), \ \left<\varphi^*\right>_{\hat\Om} = 1.
\end{equation}
We can solve (\ref{C3}), (\ref{D3}) uniquely by standard methods.  Thus putting $D\varphi^*(x) = \varphi^*(x) -\varphi^*(x-1)$, $1 \le x \le L$, then we may write (\ref{C3}) as
\begin{equation} \label{E3}
\frac 1 2 \Big[ D\varphi^*(x) - D\varphi^*(x+1)\Big] - b(x-1) \varphi^*(x-1)+b(x+1)\varphi^*(x+1)=0, \ 1 \le x \le L.
\end{equation}
If we sum (\ref{E3}) over the set $\{1 \le x \le y\}$ we obtain the equation,
\begin{multline*}
\frac 1 2 \Big[ D\varphi^*(1) - D\varphi^*(y+1)\Big] - b(0)  \varphi^*(0) - b(1)  \varphi^*(1) \\
+ b(y) \varphi^*(y)+ b(y+1)\varphi^*(y+1) = 0, \quad 1 \le y \le L.
\end{multline*}
Then, using the fact that $\varphi^*(1) = \varphi^*(0)$, $b(1) = -b(0)$, we conclude that
\[	\varphi^*(y+1) = \bar\de_y \varphi^*(y) \big/ \de_{y+1}, \quad 1 \le y \le L,	\]
whence we have
\begin{equation} \label{F3}
\varphi^*(y) = \varphi^*(1) \; \prod^{y-1}_{j=1} \ \bar\de_j / \de_{j+1} \ , \ \ 1 \le y \le L.
\end{equation}
The formula (\ref{A3}) follows from (\ref{F3}) and the normalization condition in (\ref{D3}). 
\end{proof}
\begin{lemma}  Let $\hat\Om$ be the space $\hat\Om = \{ x \in \Z : 1 \le x \le L\}$.  If $\psi_0 : \hat\Om \ra \R$ is the solution to the homogeneous equation (\ref{AW2}) with the boundary conditions  (\ref{BB2}) then $\psi_0(1)$ is given by the formula
\begin{equation} \label{G3}
2\psi_0(1) =  \prod^L_{k=1} \bar\de_k \big/ \sum^L_{r=1} \ \prod^{r-1}_{j=1} \de_j \prod^L_{j=r+1} \bar\de_j,
\end{equation}
where $\de_j,\; \bar\de_j, \; 1 \le j  \le L$, are as in (\ref{B3}).
\end{lemma}
\begin{proof}
From (\ref{AU2}, (\ref{AW2}), (\ref{BB2}), we see that $\psi_0(x)$ satisfies the equation,
\begin{equation} \label{H3}
\psi_0(x) - \frac 1 2 \Big[ \psi_0(x+1) + \psi_0(x-1)\Big] - b(x)\Big[ \psi_0(x+1) - \psi_0(x-1)\Big]=0, \ 1 \le x \le L,
\end{equation}
with the boundary conditions,
\begin{equation} \label{I3}	
\psi_0(0) = -\psi_0(1), \ \ \ \psi_0(L+1) = 1 - \psi_0(L).		\
\end{equation}
We can solve (\ref{H3}), (\ref{I3}) by standard methods.  Thus putting $D\psi_0(x) = \psi_0(x) - \psi_0(x-1)$ equation (\ref{H3}) implies
\begin{equation} \label{J3}
D\psi_0(x+1)=\de_x D\psi_0(x)/\bar\de_x,\quad 1 \le x \le L.
\end{equation}
Observing from  (\ref{I3}) that $D\psi_0(1)=2\psi_0(1)$ we see from (\ref{J3}) that
\begin{equation} \label{K3}
D\psi_0(y+1)= 2\psi_0(1) \prod^y_{j=1} \de_j / \bar\de_j, \ \ 1 \le y \le L.
\end{equation}
If we sum (\ref{K3}) we obtain a formula for $\psi_0(y)$ given by
\begin{equation} \label{L3}
\psi_0(y)= 2\psi_0(1)\left\{ 1/2+\sum^y_{r=2} \prod^{r-1}_{j=1} \de_j / \bar\de_j\right\}, \quad 1 \le y \le L.
\end{equation}
Since $D\psi_0(L+1) = 1-2\psi_0(L)$ from (\ref{I3}) it follows that if we add (\ref{K3}) to twice (\ref{L3}) when $y=L$ we obtain a formula for $\psi_0(1)$ given by
\begin{equation} \label{M3}
 2\psi_0(1) =  \prod^L_{k=1} \bar\de_k \Big/  \left\{ \prod^L_{j=1} \bar\de_j +2 \sum^L_{r=2} \prod^{r-1}_{j=1} \de_j \prod^L_{j=r}\bar\de_j + \prod^L_{j=1} \de_j \right\} .
\end{equation}
One can easily see that the denominator of the expression in (\ref{M3}) can be rewritten as in (\ref{G3}). \ \end{proof}
\begin{remark}  Observe from (\ref{A3}), (\ref{G3}) that under the reflection $b \ra -b$ the expression $\varphi^*(1)\de_1 / L$ becomes $2\psi_0(1)$.
\end{remark}
\begin{lemma}  There is the inequality $\varphi^*(1)\de_1  \psi_0(1) \le 1/8L$.
\end{lemma}
\begin{proof}
For $1 \le r \le L$ define $a_r$ by
\[	a_r = \prod^{r-1}_{j=1} \bar\de_j \ \prod^{L}_{j=r+1} \de_j ,	\]
and $\bar a_r$ the corresponding value of $a_r$ under the reflection $b \ra -b$.  From Lemma 3.1, 3.2 we see that we need to prove that 
\[	4L^2 \ \prod^L_{k=1} \de_k \bar\de_k \le \left\{  \sum^L_{r=1}  a_r \right\}\left\{ \sum^L_{r=1} \bar a_r \right\}. \]
Using the fact that for $1 \le r,s \le L$, 
\[	(a_r \bar a_s + a_s \bar a_r)/2 \ge (a_r \bar a_r a_s \bar a_s)^{1/2},	\]
we see that
\[  \left\{  \sum^L_{r=1}  a_r \right\}\left\{ \sum^L_{r=1} \bar a_r \right\} \ge \sum^L_{r,s=1} (a_r \bar a_r  a_s \bar a_s)^{1/2} \]
\[ \ge  \sum^L_{r,s=1} 4 \ \prod^L_{k=1} \de_k \bar\de_k = 4L^2 \prod^L_{k=1} \de_k \bar\de_k,   \]
where we have used the fact that for $1 \le j \le L$, one has $\de_j \bar\de_j \le 1/4$. 
\end{proof}
\begin{proof}[Proof of Theorem 1.2]
 $(d=1)$: This follows from Lemmas 3.1 - 3.3 and Lemma 2.5, using the RHS of (\ref{BC2}) as the representation for $q(b)$. 
 \end{proof} 
Next we turn to the $d > 1$ case with $L_1 = 2$.  Then we can write $\hat \Om = \{(0,y):y\in\Om_{d-1}\}$ where $\Om_{d-1} \subset \Z^{d-1}$ is a $d-1$ dimensional rectangle.  It is easy to see now from (\ref{BK2}), on using the anti-symmetry of $b$ and the symmetry of $\varphi^*$, that $\varphi^* \equiv 1 $.  Also from (\ref{BM2}), on using the anti-symmetry of $b$ and $\varphi$, we have that
\begin{equation} \label{N3}
\varphi(0,y) = 2d \left[ -\De_{d-1} +4\right]^{-1} b(0,y),
\end{equation}
where in (\ref{N3}) the operator $\De_{d-1}$ is the discrete Laplacian acting on functions $\Psi : \Om_{d-1} \ra \R$.  We have then from (\ref{W2}) that $\psi(0,y) = -2b(0,y)\varphi(0,y)$, and so we get the formula for the effective diffusion constant,
\begin{equation} \label{O3}
1/2d+ 2 <\phi^*\psi> = 1/2d - 8d \left< b(\cd)\left[-\De_{d-1} + 4\right]^{-1} b\right>_{\Om_{d-1}}.
\end{equation}
It is clear that the RHS of (\ref{O3}) is smaller than $1/2d$.  We can alternatively derive the effective diffusion constant formula by using the expression on the RHS of (\ref{BC2}).  Thus we have
\[   \psi_0(0,y) = 2d \left[ -\De_{d-1} + 4\right]^{-1} [ 1/2d +b(0,y)],   \]
whence the effective diffusion constant is given by the formula
\begin{multline} \label{P3}
L^2_1 \left< \phi^*(\cd) \left[ 1/2d - b(\cd)\right] \psi_0(\cd) \chi_0(\cd) \right>_{\hat\Om} =
\\
  8d \left< \left[ 1/2d - b(\cd)\right] \left[-\De_{d-1} + 4\right]^{-1} \left[ 1/2d + b(\cd)\right] \right>_{\Om_{d-1}}. 
\end{multline}

We shall use the formula on the LHS of (\ref{P3}) to obtain an expression for the effective diffusion constant in the case 
$L_1=4$.  Here $\hat \Om$ is the space $\hat \Om = \{(n,y) : n=0,1, \ y \in \Om_{d-1}\}$.  For  $y \in \Om_{d-1}$ we define $\de_y, \; \bar\de_y, \; \ve_y, \; \bar\ve_y$ by
\begin{eqnarray} \label{Q3}
\de_y &=&  1/2d - b(0,y), \ \bar\de_y = 1/2d + b(0,y), \\
\ve_y &=& 1/2d + b(1,y), \ \bar\ve_y = 1/2d - b(1,y). \nn 
\end{eqnarray}
We see then from (\ref{BK2}), (\ref{Q3}) that $\vp^*$ satisfies the system of equations,
\begin{eqnarray} \label{R3}
\left( \frac{-\De_{d-1} + 2}{2d}\right) \vp^*(0,y) &-& \bar\ve_y \; \vp^*(1,y) - \de_y\vp^*(0,y) = 0, \\
\left( \frac{-\De_{d-1} + 2}{2d}\right) \vp^*(1,y) &-& \ve_y \; \vp^*(1,y) - \bar\de_y\vp^*(0,y) = 0. \nn
\end{eqnarray}
Adding the 2 equations above we conclude that $-\De_{d-1}[\vp^*(0,y) + \vp^*(1,y)]=0$, $y\in\Om_{d-1}$, whence on using the normalization $<\vp^*>_{\hat\Om} = 1$ we conclude that $\vp^*(0,y) + \vp^*(1,y) = 2, \ y \in \Om_{d-1}$.  Hence from (\ref{R3}) we have that $\vp^*(0,y)$ satisfies the equation,
\begin{equation} \label{S3}
\left[ -\De_{d-1}/2d + \bar\ve_y +\bar\de_y \right] \vp^*(0,y) = 2\bar\ve_y, \ y \in \Om_{d-1}.
\end{equation}
Evidently (\ref{S3}) has a unique positive solution.

We proceed similarly to obtain a formula for $\psi_0$.  Thus $\psi_0$ satisfies the system of equations,
\begin{eqnarray} \label{T3}
\left( \frac{-\De_{d-1} + 2}{2d}\right) \psi_0(0,y) &-& \bar\de_y \;\psi_0(1,y) + \de_y \psi_0(0,y) = 0, \\
\left( \frac{-\De_{d-1} + 2}{2d}\right) \psi_0(1,y) &+& \ve_y \psi_0(1,y) - \bar\ve_y \psi_0(0,y) = \ve_y, \quad y \in \Om_{d-1}. \nn
\end{eqnarray}
Adding the two equations in (\ref{T3}) we get 
\begin{equation} \label{U3}
\left[ -\De_{d-1}/2d + \ve_y +\de_y \right] \left\{ \psi_0(0,y) + \psi_0(1,y)\right\}=\ve_y, \quad y \in \Om_{d-1}.
\end{equation}
We may also rewrite the first equation of (\ref{T3}) as,
\begin{equation} \label{V3}
\left( \frac{-\De_{d-1} + 4}{2d}\right) \psi_0(0,y) = \bar\de_y \{\psi_0(0,y) +  \psi_0(1,y)\}, \quad y \in \Om_{d-1}.
\end{equation}
We conclude from (\ref{S3}), (\ref{U3}), (\ref{V3}),  and the formula on the LHS of (\ref{P3}) that the effective diffusion constant $q(b)$ is given by,
\begin{multline} \label{W3}
q(b) = 2^7d^3 \Big< \left\{ \de \left[ -\De_{d-1} + 2 - V\right]^{-1} \bar\ve \right \}  \\ 
\left[ -\De_{d-1} + 4  \right]^{-1} 
\left\{ \bar\de \left[ -\De_{d-1} + 2 + V\right]^{-1} \ve \right \} \Big>_{\Om_{d-1}},
\end{multline}
where $V : \Om_{d-1} \ra \R$ is given by $V(y) = 2d[b(1,y) - b(0,y)], \ y \in \Om_{d-1}$.

We first show that $q(b) \le 1/2d$ in the case where $V$ is constant.
\begin{lemma}  Let $q(b)$ be given by (\ref{W3}) and assume $V$ is constant.  Then there is the inequality, $q(b) \le 1/2d.$
\end{lemma}
\begin{proof}
Since $\ve + \de= (2 + V)/2d$ there is an $f : \Om_{d-1} \ra \R$ such that 
\begin{eqnarray} \label{X3}
\ve &=& (2+V)/4d + f, \ \ \de = (2+V)/4d - f, \\
\bar \ve &=& (2-V)/4d - f, \ \ \bar\de = (2-V)/4d + f. \nn
\end{eqnarray}
We rewrite the expression on the RHS of (\ref{W3}) in terms of $f$.  To do this we let $w_+, w_-$ be solutions to the equations,
\begin{eqnarray} \label{Y3}
\left[ -\De_{d-1} + 2 + V \right]w_+ &=& f, \\
\left[ -\De_{d-1} + 2 - V \right]w_- &=& f. \nn
\end{eqnarray}
It follows that
\begin{eqnarray*}
\left[ -\De_{d-1} + 2 + V \right]^{-1} \ve &=& 1/4d + w_+ , \\
\left[ -\De_{d-1} + 2 - V \right]^{-1} \bar\ve &=& 1/4d - w_- . 
\end{eqnarray*}
Hence from (\ref{W3}) $q(b)$ is given by the expression,
\begin{multline*}
q(b) = 2^7d^3  \Bigg<  \left\{ \left[ \frac{2+V}{4d} - f\right]\left[ \frac 1{4d} - w_-\right]\right\} \\
\left[-\De_{d-1} + 4\right]^{-1} 
 \left\{ \left[ \frac{2-V}{4d} + f\right]\left[ \frac 1{4d} + w_+\right]\right\} \Bigg>_{\Om_{d-1}}.
\end{multline*}
This is a quartic expression in $f$ and the zeroth order term is given by,
\begin{equation} \label{Z3}
{\rm zeroth \ order} \ \ = \frac 1{2d} \left< (2+V)\left[ -\De_{d-1}+4\right]^{-1} (2-V) \right>_{\Om_{d-1}}.
\end{equation}
Observe that the expression in (\ref{Z3}) is identical to the RHS of (\ref{P3}) if $\ve = \de$.  For the first order term we have the expression,
\begin{multline} \label{AA3}
{\rm first \ order} \ \ = 2 \left< (2+V)\left[ -\De_{d-1}+4\right]^{-1} f \right>_{\Om_{d-1}} - 
\ 2 \left< f \left[ -\De_{d-1}+4\right]^{-1} (2-V) \right>_{\Om_{d-1}}  \\
+ 2 \left< (2+V)\left[ -\De_{d-1}+4\right]^{-1} (2-V)w_+ \right>_{\Om_{d-1}}   
-2 \left< (2+V) w_- \left[ -\De_{d-1}+4\right]^{-1} (2-V) \right>_{\Om_{d-1}} .   
\end{multline}
Now from (\ref{Y3}) we have that 
\begin{eqnarray} \label{AC3}
 (2-V)w_+ &=& \left[ -\De_{d-1}+4\right] w_+ - f, \\
 (2+V)w_- &=& \left[ -\De_{d-1}+4\right] w_- - f. \nn
\end{eqnarray}
From this we conclude that the expression in (\ref{AA3}) is the same as
\[	2 \left< (2+V)w_+ \right>_{\Om_{d-1}}  - 2 \left< (2-V)w_- \right>_{\Om_{d-1}} = 0.  \]
The second order term in (\ref{W3}) is given by
\begin{multline} \label{AB3}
- 2^3d \left< f \left[ -\De_{d-1}+4\right]^{-1} \ f \right>_{\Om_{d-1}} 
- 2^3d \left< (2+V)w_- \left[ -\De_{d-1}+4\right]^{-1} (2-V)w_+ \right>_{\Om_{d-1}}  \\ 
- 2^3d \left< f \left[ -\De_{d-1}+4\right]^{-1} (2-V)w_+ \right>_{\Om_{d-1}} 
- 2^3d \left< (2+V)w_- \left[ -\De_{d-1}+4\right]^{-1} f \right>_{\Om_{d-1}}   \\
+ 2^3d \left< (2+V) \left[ -\De_{d-1}+4\right]^{-1}  fw_+ \right>_{\Om_{d-1}}  
+ 2^3d \left< f w_- \left[ -\De_{d-1}+4\right]^{-1} (2-V) \right>_{\Om_{d-1}} .  
\end{multline}
Observe from (\ref{AC3}) that there is the identity,
\begin{multline} \label{AD3}
 \left< (2+V)w_- \left[ -\De_{d-1}+4\right]^{-1} (2-V)w_+ \right>_{\Om_{d-1}} 
=  \left< f \; \left[ -\De_{d-1}+4\right]^{-1} \; f \right>_{\Om_{d-1}}  \\
+ \left< w_-\left[ -\De_{d-1}+4\right]w_+ \right>_{\Om_{d-1}}  
-  \left< f w_+ \right>_{\Om_{d-1}} - \left< f w_- \right>_{\Om_{d-1}} . 
\end{multline}

We similarly have that 
\begin{equation} \label{AE3}
 \left< f \; \left[ -\De_{d-1}+4\right]^{-1} (2-V)w_+ \right>_{\Om_{d-1}} = 
-  \left< f \; \left[ -\De_{d-1}+4\right]^{-1} \; f \right>_{\Om_{d-1}} + \left< f w_+ \right>_{\Om_{d-1}} ,  
\end{equation}
$$
\left< (2+V)w_- \left[ -\De_{d-1}+4\right]^{-1} \; f \right>_{\Om_{d-1}} =
-  \left< f \; \left[ -\De_{d-1}+4\right]^{-1} \; f \right>_{\Om_{d-1}} + \left< f w_- \right>_{\Om_{d-1}} . 
$$
Define now $U : \Om_{d-1} \ra \R$ as the solution to the equation,
\begin{equation} \label{AF3}
\left[ -\De_{d-1}+4\right]U = V.
\end{equation}
Then from equations (\ref{AB3}) - (\ref{AF3}) we see that the second order term in (\ref{W3}) is given by
\begin{eqnarray} \label{AG3}
{\rm second \ order} &=& 2^3d \Bigg\{ - \left< w_- \left[ -\De_{d-1}+4\right]w_+ \right>_{\Om_{d-1}} \\
&+& \frac 1 2 \left< f \left[ w_+ + w_-\right]\right>_{\Om_{d-1}} + \left< fU \left[ w_+ - w_-\right]\right>_{\Om_{d-1}} \Bigg\} . \nn
\end{eqnarray}
The term of third order is given by
\begin{multline} \label{AH3}
2^5d^2 \left< f\; w_- \left[ -\De_{d-1}+4\right]^{-1} \ f \right>_{\Om_{d-1}} + 
2^5d^2 \left< f\; w_- \left[ -\De_{d-1}+4\right]^{-1} (2-V)w_+ \right>_{\Om_{d-1}} \\
-  2^5d^2 \left< f \left[ -\De_{d-1}+4\right]^{-1}f\; w_+\right> - 2^5d^2 \left< (2+V)w_- \left[ -\De_{d-1}+4\right]^{-1} \; fw_+ \right>_{\Om_{d-1}}. 
\end{multline}
Using (\ref{AC3}) again we see that the expression (\ref{AH3}) is the same as
\[   2^5d^2 \left< f\; w_- \; w_+ \right>_{\Om_{d-1}} - 2^5d^2 \left< f\; w_- \; w_+ \right>_{\Om_{d-1}} = 0. \]
Finally the fourth order term is given by
\begin{equation} \label{AI3}
{\rm fourth \ order} \ = 2^7d^3 \left< f\; w_- \left[ -\De_{d-1}+4\right]^{-1} f\; w_+ \right>_{\Om_{d-1}} .
\end{equation}
Hence from (\ref{Z3}), (\ref{AG3}), (\ref{AI3}) we have the formula for $q(b)$,
\begin{multline} \label{AJ3}
q(b) =  \frac 1{2d} \left< (2+V)\left[ -\De_{d-1}+4\right]^{-1} (2-V) \right>_{\Om_{d-1}} 
- 2^3d \left< w_- \left[ -\De_{d-1}+4\right]w_+ \right>_{\Om_{d-1}} \\
+ 2^2d \left< f[w_+ + w_-]\right>_{\Om_{d-1}} 
+ 2^3d \left< fU [w_+ - w_-]\right>_{\Om_{d-1}} + 2^7d^3 \left< f\; w_- \left[ -\De_{d-1}+4\right]^{-1} f\; w_+ \right>_{\Om_{d-1}}. 
\end{multline}
It is clear from the definitions of $V,f$ that
\begin{equation} \label{AK3}
|V(y)| < 2, \ |f(y)| < [2 - |V(y)|] /4d, \ \ \ y \in \Om_{d-1}.
\end{equation}
From (\ref{AK3}) it follows that there is the inequality,
\begin{equation} \label{AL3}
2^7d^3 \left< fw_- \left[ -\De_{d-1} + 4\right]^{-1} fw_+\right>_{\Om_{d-1}} \le 2^4d^3 \left< (fw_-)^2\right>_{\Om_{d-1}} 
\end{equation}
\[  + 2^4d^3 \left< (fw_+)^2\right>_{\Om_{d-1}} \le d \left< [2 - |V|]^2 [w^2_- + w^2_+] \right> _{\Om_{d-1}}. \]
We define now a quadratic form $Q_V$ depending  on $V$ by
\begin{multline} \label{AM3}
Q_V(f) = \left< w_-\left[ -\De_{d-1} + 4\right] w_+ \right>_{\Om_{d-1}} 
- \frac 1 2 \left< f [w_+ + w_-] \right> _{\Om_{d-1}}  \\
- \left< fU [w_+ - w_-] \right> _{\Om_{d-1}} 
- \frac 1 8 \left< [2 - |V|]^2 [w^2_- + w^2_+] \right> _{\Om_{d-1}}. 
\end{multline}
It is evident from (\ref{AJ3}), (\ref{AL3}) that
\begin{equation} \label{AN3}
q(b) \le \frac 1{2d} \left< (2+V)\left[ -\De_{d-1} + 4\right]^{-1} (2-V)\right>_{\Om_{d-1}}- 2^3 d \ Q_V(f).
\end{equation}
Thus to prove the result it will be sufficient to show that $Q_V$ is nonnegative definite.  

Since $V$ is constant we can compute the eigenvalues of $Q_V$ explicitly.  Thus if $p^2$ denotes the eigenvalue of $-\De_{d-1}$, corresponding to the eigenvector $\exp[i\xi \cd x], x \in \Om_{d-1}$, then the eigenvalue of $Q_V$ is
\begin{multline} \label{AO3}
\bigg\{ 2(4 + V^2)(p^2 + 2 + V)(p^2 + 2 - V) -
\\
	\frac 12(2 - |V|)^2 \left[ (p^2 + 2 + V)^2+(p^2 + 2 - V)^2\right] \bigg\} \bigg/ 4(p^2 + 2 + V)^2 (p^2 + 2 - V)^2.  
\end{multline}
We can rewrite the numerator of (\ref{AO3}) as
\begin{equation} \label{AP3}
\left\{ 2(4 + V^2) - (2 - |V|)^2 \right\} [p^2+2]^2- V^2 \left\{ 2(4 + V^2) + (2 - |V|)^2 \right\}.
\end{equation}
Since $|V| < 2$ the expression in (\ref{AP3}) is bounded below by its value for $p=0$ which can be written as
\[	(4 + V^2)  (2 - |V|) [2 + 3 |V|] \ge 0.	\]
\end{proof}

We proceed now to the general case which will follow from the fact that the quadratic form (\ref{AM3}) is nonnegative definite for any $V$ satisfying (\ref{AK3}). From here on we shall denote $\De_{d-1}, \left< \cd\right >_{\Om_{d-1}}$ simply as 
$\De$ and $\left< \cd \right>$ respectively.  We first note that by using (\ref{Y3}) we can obtain some alternate representations for $Q_V$.  Thus if we write
\[ \left< w_-\left[ -\De_{d-1} + 4\right] w_+ \right> = 2 \left<w_- \; w_+ \right> + \frac 1 2 \left< f[w_+ + w_-] \right>, \]
we see that $Q_V$ is given by
\begin{equation} \label{AQ3}
Q_V(f) = 2 \left< w_- \; w_+ \right> - \left< fU [w_+ - w_-] \right> - \frac 1 8 \left< [2 - |V|]^2 [w^2_- + w^2_+]\right> .
\end{equation}
We also have that
 \begin{multline*}
	\left< f U [w_+ - w_-] \right> = \left< \left\{ U[-\De + 2 - V]w_- \right\} w_+ \right>  
 - \left< \left\{ U[-\De + 2 +V]w_+ \right\} w_- \right> \\
  = - 2 \left< UV w_- \; w_+ \right> - \left< U w_+ \; \De w_- \right> + \left< U w_- \; \De w_+ \right>.
  \end{multline*}

Hence from (\ref{AQ3}) we have the formula,
\begin{equation} \label{AR3}
Q_V(f) = 2 \left< w_- \; w_+[1+UV] \right> + \left< U w_+ \; \De w_- \right> -\left< U w_- \; \De w_+ \right> - \frac 1 8 \left< [2 - |V|]^2 [w^2_- + w^2_+]\right> .
\end{equation}

We first show that a simple quadratic form related to $Q_V$ is nonnegative definite.

\begin{lemma} Let $V$ satisfy (\ref{AK3}) and $w_+ \;, w_-$ be solutions to (\ref{Y3}) for any $f : \Om_{d-1} \ra \R$.  Then there is the inequality $\left<w_+ \; w_- \right>  \ \ge 0$.
\end{lemma}
\begin{proof}
Let $w$ be the solution to the equation,
\begin{equation}  \label{AS3}
[-\De + 2 + V] [-\De + 2 - V]w = f.
\end{equation}
 Then from (\ref{Y3}) we see that $w_+ = [-\De + 2 - V]w$.  Hence we have from (\ref{Y3}), (\ref{AS3}) that
\begin{multline*}
\left< w_+ \; w_-\right> = \left< [-\De + 2 - V]w \ w_- \right> = \left< w[-\De + 2 - V] w_- \right> \\
 = \left<w \ f \right> =
\left<\left\{ \left[ -\De + 2 + V \right] w \right\} \left\{ \left[-\De + 2 - V \right] w \right\} \right> \\
=\left< \left\{ (-\De + 2 )w \right\}^2 \right> - \big< V^2w^2 \big> \ge \big< (4-V^2)w^2 \big> \ge 0.
\end{multline*}
\end{proof}
To proceed further we need to localize the quadratic form (\ref{AR3}).
\begin{lemma} Let  ${\cL}_+, {\cL}_-$ be operators on functions $\Phi : \Om_{d-1} \ra \R$ defined by 
\begin{eqnarray*}
{\cL}_+ \Phi &=& (-\De + 2) \Phi / V - \Phi, \\
{\cL}_- \Phi &=& (-\De + 2) \Phi / V + \Phi, 
\end{eqnarray*}
where we assume $V$ satisfies (\ref{AK3}) and $V(y) \not= 0, y \in \Om_{d-1}$.  Then ${\cL}_+, {\cL}_-$ are invertible and there is the identity,
\begin{equation} \label{AT3}
[-\De + 2 + V]{\cL}_+=    [-\De + 2 - V] {\cL}_- \ .
\end{equation}
\end{lemma}
\begin{proof}
Verification. 
\end{proof}
It follows from Lemma 3.6 that we can choose $f$ in (\ref{Y3}), (\ref{AR3}) as the operator (\ref{AT3}) acting on a function $\Phi : \Om_{d-1} \ra \R$, in which case $w_+ = {\cL_+
}\Phi, \ w_- = {\cL}_- \Phi$.  If we substitute into 
(\ref{AR3}) we obtain a quadratic form  $\tilde Q_V(\Phi)$ which is local in $\Phi$, and it is this quadratic form which we will show is nonnegative definite.  First we show that the quadratic form obtained from $\tilde Q_V$ upon replacing $U$ by $V/4$ is nonnegative definite.
\begin{lemma} For $\Phi : \Om_{d-1} \ra \R$ and $w_+ = {\cL}_+ \Phi$, $w_- = {\cL}_- \Phi$, there is for $d=2$ the inequality,
\begin{eqnarray} \label{AU3}
2 \left< w_- \; w_+ [4 + V^2] \right> &+& \left< Vw_+ \; \De w_- \right> - \left< Vw_- \; \De w_+ \right> \\
&-& \frac 1 2  \left< [2 - |V|]^2 [w^2_- + w^2_+]\right> \ge 0. \nn
\end{eqnarray}
\end{lemma}
\begin{proof}
We first note that the first term in (\ref{AU3}) is nonnegative.  Thus we have
\begin{multline*}
\left< w_- \; w_+ [4 + V^2] \right> = \left< [(-\De + 2)\Phi]^2 [4/V^2 + 1] \right> 
- \left<\Phi^2 [4 + V^2] \right>  \\
\ge 2 \left< [(-\De + 2)\Phi]^2 \right> - 
 \left< \Phi^2[4 + V^2] \right> \ge  \left< \Phi^2[4 - V^2] \right> \ge 0,
\end{multline*}
where we have used (\ref{AK3}).  The second and third terms of (\ref{AU3}) are given by the formula,
\begin{multline*}
\left< Vw_+ \; \De w_- \right> - \left< Vw_- \; \De w_+ \right> =  \\
2\left< [(-\De + 2)\Phi](\De \; \Phi) \right> - 2\left< [\De(V\Phi)] \left[ \frac{(-\De +2)\Phi}{V} \right]\right>, 
\end{multline*}
on summation by parts.  From the last two equations we therefore have that
\begin{multline} \label{AV3}
 2\left< w_- \; w_+ [4 + V^2] \right> + \left< Vw_+ \; \De w_- \right> - \left< Vw_- \; \De w_+ \right> \\
= 8 \left< \frac{[ (-\De +2)\Phi]^2} {V^2} \right> - 2 \left< V^2 \Phi^2\right> + 4 \left< (\nabla \Phi)^2 \right> 
- 2 \left< [ \De(V \Phi)] \left[ \frac {(-\De +2)\Phi]}{V} \right] \right> . 
\end{multline}
Using the fact that
\begin{multline*}
-\frac{1} {V(x)} \ \De(V \Phi)(x) = 2(d-1) \Phi(x) \\
- \sum^d_{j=2} \left[ \frac{V(x+{\bf e}_j)} {V(x)} \Phi(x+{\bf e}_j) + \frac{V(x-{\bf e}_j)}{V(x)} \Phi(x-{\bf e}_j)\right],
\end{multline*}
we conclude from (\ref{AV3}) that
\begin{multline} \label{AW3}
\left< w_- \; w_+ [4 + V^2] \right> + \frac 1 2 \left< Vw_+ \; \De w_- \right> - \frac 1 2 \left< Vw_- \; \De w_+ \right> \\
= 4 \left< \frac{[ (-\De +2)\Phi]^2}{V^2} \right> + \left< [4(d-1) - V^2]\Phi^2 \right> + 2d \left< (\nabla \Phi)^2 \right>  \\
- \left< [(-\De + 2)\Phi(\cd)] \sum^d_{j=2} \left[ \frac{V(\cd +{\bf e}_j)}{V(\cd)} \Phi(\cd +{\bf e}_j) + \frac{V(\cd-{\bf e}_j)}{V(\cd)} \Phi(\cd -{\bf e}_j) \right] \right>.
\end{multline}
Now the Schwarz inequality yields
\begin{eqnarray*}
\left| \left[ (-\De+2)\Phi(x)\right] \frac{V(y)\Phi(y)}{V(x)} \right| &\le& V(y)^2 \left[ \alpha \; \frac{\Phi(y)^2}{V(y)^2} + \frac 1{4\alpha} \frac{[(-\De + 2)\Phi(x)]^2}{V(x)^2} \right] \\
&\le& \alpha \Phi(y)^2 + \frac 1\alpha \frac{[(-\De + 2)\Phi(x)]^2}{V(x)^2}, \ \ \ x,y \in \Om_{d-1},
\end{eqnarray*}
for any $\al > 0$.  Hence there is from (\ref{AW3}) the inequality,
\begin{multline} \label{AVV3}
\left< w_- \; w_+ [4 + V^2] \right> + \frac 1 2 \left< Vw_+ \; \De w_- \right> - \frac 1 2 \left< Vw_- \; \De w_+ \right> \\
\ge \left[ 4 - \frac{2(d-1)}{\al} \right] \left< \frac{[ (-\De +2)\Phi]^2}{V^2} \right> +  \\
 \left< \left[2(d-1)(2-\al)-V^2\right]\Phi^2\right> + 2d\left<(\nabla\Phi)^2\right>. 
\end{multline}
For $d=3$ and $\al = 1$ the RHS of (\ref{AVV3}) is evidently nonnegative but this is no longer the case when $d>3$.  For $\al = 1, d=2$ the RHS of (\ref{AVV3}) is bounded below by the nonnegative quantity,
\begin{equation} \label{AWW3}
\frac 1 2 \left< \left( \frac{4-V^2}{V^2}\right) \left[ (-\De + 2)\Phi\right]^2\right> + \left< (4-V^2)\Phi^2 \right>.
\end{equation}
We consider now the last term in the expression (\ref{AU3}).  We have that
\begin{equation} \label{AX3}
\frac 1 4 \left< [2 - |V|]^2 [w^2_- + w^2_+] \right> = \frac 1 2 \left< [2-|V|]^2 \left\{ \frac{[ (-\De +2)\Phi]^2}{V^2} + \Phi^2\right\} \right>.
\end{equation}
If we  now use the inequality $[2 - |V|]^2 \le 4 - V^2$ we see that the expression (\ref{AX3}) is less than (\ref{AWW3}).  Hence the inequality (\ref{AU3}) is established for $d=2$. 
\end{proof}
Next we turn to showing that $\tilde Q_V$ is nonnegative definite for $d=2$.  To do this we write the solution $U$ of (\ref{AF3}) as
\begin{equation} \label{AZ3}
U(x) = \frac 1 4 \sum_y G(y) V(x+y),
\end{equation}
where $G(y)$ is the Green's function for $[-\De/4 + 1]^{-1}$, whence $G(y)$ is nonnegative for all $y$ and
\begin{equation} \label{AY3}
\sum_y G(y) = 1, \ \ \ \ G(y) = G(-y), \ \ y = 1,2,... \ \ .
\end{equation}
We consider the first three terms in the expression (\ref{AR3}) for $Q_V$.  Using Lemma 3.6 and setting $w_+ = {\cL}_+ \Phi, w_- = {\cL}_- \Phi$ we have that
\begin{eqnarray} \label{BA3}
&\ & \left< w_-  w_+ [1 + UV] \right> + \frac 1 2 \left< Uw_+  \De w_- \right> -\frac 1 2 \left< Uw_-  \De w_+\right> \\
&=& \left< \frac 1{V^2} \left[ (-\De + 2) \Phi \right]^2 [1 + UV] \right> - \left< \Phi^2[1 + UV] \right> \nn \\
&+& \left< \frac U V [\De\Phi] \left[ (-\De + 2) \Phi\right]\right> - \left< \De(U\Phi) \frac 1 V [-\De+2]\Phi \right> \nn \\
&=& \left< \frac 1{V^2} \left[ (-\De + 2) \Phi \right]^2\right>  - \left< \Phi^2 \right> - \left<\Phi^2 UV\right> \nn \\
&+&2 \left< \frac U V \Phi (-\De + 2) \Phi \right> - \left< \De(U\Phi) \frac 1 V[-\De + 2]\Phi \right> \nn \\
&=& \left< \frac 1{V^2} \left[ (-\De + 2) \Phi \right]^2 \right> - \left< \Phi^2 \right> - \left<\Phi^2 UV\right> \nn \\
&+& 4 \left< \frac U V \Phi (-\De + 2) \Phi \right> - \left< (\tau_1 U)(\tau_1 \Phi) \frac 1 V[-\De + 2]\Phi\right> \nn \\
&-&  \left< (\tau_{-1} U)(\tau_{-1} \Phi) \frac 1 V[-\De + 2]\Phi\right>, \nn
\end{eqnarray}
where $\tau_x \vp(y) = \vp(x + y), \ y \in \Z$.  We consider the last three terms in the previous expression.  We write using (\ref{AZ3}),
\begin{multline*}
4 \left< \frac U V \Phi (-\De + 2) \Phi \right> = G(0) \left< (\nabla \Phi)^2 +2\Phi^2 \right> 
+ \sum_{y \not= 0} G(y) \left< \frac{\tau_y V}{V} \ \Phi(-\De + 2) \Phi \right>, \\
\left< (\tau_1 U)(\tau_1 \Phi) \frac 1 V[-\De + 2]\Phi\right> 
= \frac 1 4 G(-1) \left\{ \left< \nabla(\tau_1\Phi) \nabla\Phi \right> 
 + 2 \left<( \tau_1\Phi) \Phi \right> \right\} \\
+ \frac 1 4 \sum_{y \not= 0} G(y-1) \left< \frac{\tau_y V}{V} \ \tau_1 \Phi(-\De + 2) \Phi \right> ,
\end{multline*}
with a similar expression for the last term in (\ref{BA3}).  We conclude from this that the last three terms of (\ref{BA3}) are given by,
\begin{multline} \label{BB3}
4 \left< \frac U V \Phi (-\De + 2) \Phi \right> - \left< (\tau_1 U)(\tau_1 \Phi) \frac 1 V[-\De + 2]\Phi\right> 
-  \left< (\tau_{-1} U)(\tau_{-1} \Phi) \frac 1 V[-\De + 2]\Phi\right>  \\
= G(0) \left< (\nabla \Phi)^2 + 2\Phi^2 \right>  
- \frac 1 4 G(-1) \left\{ \left< \nabla(\tau_1\Phi)\nabla \Phi \right> 
+ 2 \left< (\tau_1\Phi )\Phi \right> \right\} 
- \frac 1 4 G(1) \left\{ \left< \nabla(\tau_{-1}\Phi)\nabla \Phi \right> + 2 \left< (\tau_{-1}\Phi) \Phi \right> \right\}  \\
+ \sum_{y \ge 1} \left[ G(y) - \frac 1 4 \; G(y-1) - \frac 1 4 \; G(y+1) \right] 
 \left< \frac{\tau_y V}{V} \ \tau_1 \Phi(-\De + 2) \Phi \right>  \\
 + \sum_{y \le -1} \left[ G(y) - \frac 1 4 \; G(y-1) - \frac 1 4 \; G(y+1) \right] 
 \left< \frac{\tau_y V}{V} \ \tau_{-1} \Phi(-\De + 2) \Phi \right> \\
  + \sum_{y \ge 1} G(y) \left< \frac{\tau_y V}{V} \ \left[ \Phi - \tau_1 \Phi\right] (-\De + 2) \Phi \right> 
- \frac 1 4  \sum_{y \ge 1} G(y+1) \left< \frac{\tau_y V}{V} \ \left[ \tau_{-1}\Phi - \tau_1 \Phi\right] (-\De + 2) \Phi \right>  \\
+  \sum_{y \le -1} G(y) \left< \frac{\tau_y V}{V} \ \left[ \Phi - \tau_{-1} \Phi\right] (-\De + 2) \Phi \right> 
- \frac 1 4  \sum_{y \le -1} G(y-1) \left< \frac{\tau_y V}{V} \ \left[ \tau_1\Phi - \tau_{-1} \Phi\right] (-\De + 2) \Phi \right>. 
\end{multline}
We shall use the representation (\ref{BB3}) to show that the quadratic form $Q_V$ is non-negative definite.
\begin{lemma}  Suppose the function $G(y)$ of (\ref{AZ3}) is decreasing, non-negative for $y \ge 1$, satisfies (\ref{AY3}) and the inequalities,
\begin{equation} \label{CB3}
(-\De + 2) G(y) \le 0, \ y \ge 1; \ 1 - G(0) - 2G(1) < G(1)/2; \ G(2) < G(1)/5.
\end{equation}
Then the quadratic form $Q_V$ of (\ref{AR3}) is nonnegative definite.
\end{lemma}
\begin{proof}
We estimate the terms in (\ref{BB3}) by applying the Schwarz inequality.  Before doing this we make one further simplification of terms in (\ref{BB3}).  We write
\begin{equation} \label{BC3}
 \sum_{y \ge 1} G(y) \left< \frac{\tau_y V}{V} \ \left[ \Phi - \tau_1 \Phi\right] (-\De + 2) \Phi \right> 
+ \sum_{y \le -1} G(y) \left< \frac{\tau_y V}{V} \ \left[ \Phi - \tau_{-1} \Phi\right] (-\De + 2) \Phi \right> 
\end{equation}
$$
 =  \sum_{y \ge 2} G(y) \left< \frac{\tau_y V}{V} \ \left[ \Phi - \tau_1 \Phi\right] (-\De + 2) \Phi \right> 
+ \sum_{y \le -2} G(y) \left< \frac{\tau_y V}{V} \ \left[ \Phi - \tau_{-1} \Phi\right] (-\De + 2) \Phi \right> 
$$
$$
+ G(1) \left< \tau_1 V \left[ \frac 1 V - \frac V 4\right] \left[ \Phi - \tau_1 \Phi\right] (-\De + 2) \Phi \right> 
+G(-1) \left< \tau_{-1} V \left[ \frac 1 V - \frac V 4\right] \left[ \Phi - \tau_{-1} \Phi\right] (-\De + 2) \Phi \right> 
$$
$$
- \frac{G(1)}4 \left< (\tau_1 V) V[\Phi - \tau_1 \Phi] \De \Phi \right> 
- \frac{G(-1)}4 \left< (\tau_{-1} V) V[\Phi - \tau_{-1} \Phi] \De \Phi \right> 
+ \frac{G(1)} 2 \left< (\tau_1 V)V(\nabla \Phi)^2 \right>, 
$$
where we have used the fact that $G(1) = G(-1)$.  We also rewrite the first two terms on the RHS of (\ref{BC3}) as
\begin{equation} \label{BD3}
\sum_{y \ge 2} G(y) \left< \frac{(\tau_y V-V)}{V} \ \left[ \Phi - \tau_1 \Phi\right] (-\De + 2) \Phi \right> 
 \end{equation}
 $$
 +   \sum_{y \le -2} G(y) \left< \frac{(\tau_y V-V)}{V} \ \left[ \Phi - \tau_{-1} \Phi\right] (-\De + 2) \Phi \right> 
+\sum_{y \ge 2} \ G(y) \left< (\De\Phi)^2 + 2 (\nabla \Phi)^2 \right> .
$$
We similarly rewrite the sum of the last and third last terms of (\ref{BB3}) as
\begin{equation} \label{BE3}
- \frac 1 4 \sum_{y \ge 1} G(y+1) \left< \frac{(\tau_y V-V)}{V} \ \left[ \tau_{-1}\Phi - \tau_1 \Phi\right] (-\De + 2) \Phi \right> 
\end{equation}
$$
- \frac 1 4 \sum_{y \le -1} G(y-1) \left< \frac{(\tau_y V-V)}{V} \ \left[  \tau_1\Phi - \tau_{-1} \Phi\right] (-\De + 2) \Phi \right> . 
$$
Consider now the first three terms on the RHS of (\ref{BA3}).  These can be written as
\begin{equation} \label{BF3}
\left< \left( \frac 1{V^2} - \frac 1 4 \right) [(-\De + 2) \Phi]^2 \right> + \frac 1 4 \left< \left(\De \Phi\right)^2 \right> + \left<\left(\nabla \Phi\right)^2 \right> - 
\end{equation}
$$
\frac 1 4 \ G(0) \left< V^2 \Phi^2 \right> 
+ \frac 1 8 \sum_{y\not= 0} G(y) \left[ \left< (\tau_y V-V)^2 \Phi^2\right> - \left< \left\{ (\tau_y V)^2 + V^2 \right\}\Phi^2 \right> \right]. \nn
 $$

Next we apply the Schwarz inequality to terms in (\ref{BB3}).  Thus we estimate
\begin{equation} \label{BG3}
\left| \left< \frac{\tau_y V}{V} \; \tau_1\Phi(-\De + 2)\Phi \right> \right| \le \ \left <\Phi^2\right> + \left< \frac 1{V^2} [(-\De +2)\Phi]^2\right>,
\end{equation}
with a similar estimate when $\tau_1 \Phi$ is replace by $\tau_{-1} \Phi$.

Observe now that 
\[
\sum_{y \ge 1} \left[ G(y) - \frac 1 4 \; G(y-1) - \frac 1 4 \; G(y+1)\right] = -\frac{1}4 [2G(0)-1-G(1)],
\]
where we have used (\ref{AY3}).  Hence on using the fact that $[-\De + 2]G(y) < 0$, $y \ge 1$, we see from (\ref{BG3}) that the sum of the first five terms on the RHS of (\ref{BB3}) are bounded below by the expression,
\begin{equation} \label{BH3}
[G(0) - \frac 1 2 \; G(1)] \left< (\na \Phi)^2 + 2\Phi^2 \right> 
\end{equation}
\[	-\left\{ \left<\Phi^2\right> + \left< \frac 1{V^2} [(-\De +2)\Phi]^2\right> \right\} [G(0) - \frac 1 2 - \frac 1 2 \; G(1)].   \]
If we combine the estimate (\ref{BH3}) with (\ref{BF3}) and use the fact that $|V| < 2$ we get a lower bound for the sum of the first three terms of (\ref{BA3}) and the first five terms of (\ref{BB3}).  It is given by,
\begin{equation}\label{BI3}
\left< \left( \frac 1{V^2} - \frac 1 4 \right) \left[ (-\De + 2)\Phi\right]^2 \right> \left\{ \frac 3 2 + \frac 1 2 \; G(1) - G(0) \right\} 
\end{equation}
$$
+ \frac 1 4 \left< (\De \Phi)^2 \right> \left\{ \frac 3 2 + \frac 1 2 \; G(1) - G(0) \right\} 
+ \frac 3 2 \left< (\na \Phi)^2 \right> + \frac 1 4 \ G(0) \left<(4-V^2) \Phi^2 \right> 
$$
$$
+ \frac 1 8 \sum_{y\not= 0} G(y) \left< (\tau_y V-V)^2 \Phi^2 \right> + [1-G(0)] <\Phi^2 > 
-  \frac 1 8 \sum_{y\not= 0} G(y) \left< \left\{ (\tau_y V)^2 + V^2\right\} \Phi^2 \right>. 
$$
Observe that all terms in (\ref{BI3}) except for the final one are nonnegative.  Furthermore, the sum of the last two terms is nonnegative.

Next we estimate the terms on the RHS of (\ref{BC3}) which involve $G(1)$ and $G(-1)$.  To do this we use the Schwarz inequalities,
\begin{eqnarray} \label{BJ3}
&\ & \left| \left< \tau_1 V \left[ \frac 1 V - \frac V 4 \right] \left[ \Phi - \tau_1 \Phi \right] (-\De + 2)\Phi \right> \right| \\
&\le& \al_1 \left< \left[ \frac 1 {V^2} - \frac 1 4 \right] \left[  (-\De + 2)\Phi\right]^2\right> + \frac 1{\al_1}  
\left< (\na\Phi)^2\right>, \nn \\
&\frac 1 2 & \left| \left< (\tau_1 V)V \left[ \Phi - \tau_1 \Phi \right]  \De\Phi \right> \right| \le \al_2 \left< (\De\Phi)^2\right> +\frac 1{\al_2} \left< (\na\Phi)^2\right>\nn , 
\end{eqnarray}
for any constants $\al_1, \al_2 > 0$.  Hence on using the fact that $|V| < 2$ we see that the expression (\ref{BI3}) plus the terms in $G(1), G(-1)$ of (\ref{BC3}) is bounded below by the expression,
\begin{equation} \label{BK3}
 \left< \left( \frac 1{V^2} - \frac 1 4 \right) \left[ (-\De + 2)\Phi\right]^2 \right> \left\{ \frac 3 2 + \frac 1 2 \; G(1) - G(0) - 2\al_1 G(1) \right\} 
 \end{equation}
 $$
+ \frac 1 4 \left< (\De \Phi)^2 \right> \left\{ \frac 3 2 + \frac 1 2 \; G(1) - G(0) - 4\al_2 G(1)\right\} 
+ \left< (\na \Phi)^2 \right>  \left\{ \frac 3 2 - 2G(1) - \frac 2{\al_1} G(1) - \frac 1{\al_2} G(1) \right\} 
$$
$$
+ \frac 1 4 \ G(0) \left<(4-V^2) \Phi^2 \right> + \frac 1 8 \sum_{y\not= 0} G(y) \left<  (\tau_y V-V)^2  \Phi^2 \right>. 
$$

Let us assume for the moment that $G(y) = 0$ for $y \ge 2$.  Then (\ref{BA3}) is bounded below by (\ref{BK3}).  We write $G(0) = 1 - \ga$ whence $G(1) = \ga/2$.  Since $(-\De + 2) G(y) \le 0$, $y \ge 1$, we must have $\ga < 1/3$.  We choose $\al_1$ such that
\begin{equation} \label{BL3}
 \frac 3 2 + \frac 1 2 \; G(1) - G(0) - 2\al_1 G(1) = \frac 1 2,
\end{equation}
which yields $\al_1 = 5/4$.  We choose $\al_2$ so that
\begin{equation} \label{BM3}
\frac 3 2 + \frac 1 2 \; G(1) - G(0) - 4\al_2 G(1) = 0,
\end{equation}
which yields $\al_2 = 5/8 + 1/4\ga$.  The coefficient of $<(\na \Phi)^2>$ in (\ref{BK3}) is therefore bounded below by  $1.5 - 2.6\ga > 0$ since $\ga < 1/3$.  Hence from (\ref{AX3}) the quadratic form 
$Q_V/2$ of (\ref{AR3}) is bounded below by twice the expression, 
\begin{equation} \label{BN3}
 \frac 1 2  \left< \left( \frac 1{V^2} - \frac 1 4 \right) \left[ (-\De + 2)\Phi \right]^2 \right> + \frac 1 4 \; G(0) \left< (4-V^2)\Phi^2 \right> 
 \end{equation}
 $$
-  \frac 1 8  \left< [2 - |V|]^2 \left\{ \frac {[(-\De+2)\Phi]^2} {V^2} + \Phi^2 \right\} \right>. 
$$
If we now use the fact that $ [2 - |V|]^2 \le 4 - V^2$ we see that (\ref{BN3}) is nonnegative.

To complete the proof of the lemma we need to estimate the sum of the terms in (\ref{BD3}), (\ref{BE3}).  We rewrite these as
\begin{eqnarray} \label{BO3}
&-& \frac 1 4\; G(2) \left< \frac{(\tau_1 V-V)}{V} \ \left[ \tau_{-1}\Phi - \tau_1 \Phi\right] (-\De + 2) \Phi \right> \\
&-& \frac 1 4  \; G(-2) \left< \frac{(\tau_{-1} V-V)}{V} \ \left[  \tau_1\Phi - \tau_{-1} \Phi\right] (-\De + 2) \Phi \right>  \nn \\
&+& \sum_{y \ge 2} \left[ G(y) - \frac 1 4 \;G(y+1)\right] \left< \frac{(\tau_y V-V)}{V} \ \left[ \Phi -\tau_1\Phi \right] (-\De + 2) \Phi \right> \nn \\
&+& \frac 1 4 \sum_{y \ge 2}G(y+1)\left< \frac{(\tau_y V-V)}{V}\ \left[ \Phi - \tau_{-1} \Phi\right] (-\De + 2) \Phi \right>\nn \\
&+& \sum_{y \le -2} \left[ G(y) - \frac 1 4 \;G(y-1)\right] \left< \left(\frac{\tau_y V-V}{V}\right) \left[ \Phi -\tau_{-1}\Phi \right] (-\De + 2) \Phi \right> \nn \\
&+& \frac 1 4 \sum_{y \le -2}G(y-1)\left< \frac{(\tau_y V-V)}{V} \ \left[ \Phi - \tau_1\Phi\right] (-\De + 2) \Phi \right> \nn \\
&+& \frac 1 2 [1 - G(0) - 2G(1)]  \left< (\De \Phi)^2 + 2(\na \Phi)^2 \right> . \nn
\end{eqnarray} 
We first estimate the third term in (\ref{BO3}).  Thus we write
\begin{equation} \label{BP3}
 \left< \frac{(\tau_y V-V)}{V} \ \left[ \Phi -\tau_1\Phi \right] (-\De + 2) \Phi \right> 
 \end{equation}
 $$
= \left< (\tau_y \; V-V) \left[ \frac 1 V - \frac V 4 \right] [\Phi - \tau_1\Phi] (-\De + 2) \Phi \right> 
$$
$$
- \frac 1 4 \left< (\tau_y \; V-V) V[\Phi - \tau_1\Phi] \De  \Phi \right> 
+ \frac 1 2 \left< (\tau_y \; V-V) V[\Phi - \tau_1\Phi]   \Phi \right> . 
$$

We estimate the first two terms on the RHS of (\ref{BP3}) similarly to (\ref{BJ3}).  For the third term we use
\begin{equation} \label{BQ3}
| \left< (\tau_y \; V-V) V[\Phi - \tau_1\Phi]  \Phi \right> | \le \al \left< (\na \Phi)^2 \right> + \frac 1 \al 
\left< (\tau_y \; V-V)^2 \Phi^2 \right>, 
\end{equation}
for any $\al > 0$.  Choosing $\al = 4$ in (\ref{BQ3}) it follows that the sum of the third, fourth, fifth and sixth terms of (\ref{BO3}) is bounded below by
\begin{eqnarray} \label{BR3}
&-& \sum_{|y| \ge 2} G(y) \bigg\{ 2\al_3 \left< \left( \frac 1{V^2} - \frac 1 4 \right) [(-\De + 2)\Phi]^2 \right> \\
&+& \al_4 \left< (\De\Phi)^2 \right> + \left[ \frac 2{\al_3} + \frac 1{\al_4} + 2 \right] \left< (\na \Phi)^2 \right> + \frac 1 8 
\left< (\tau_y \; V-V)^2 \Phi^2 \right> \bigg\}, \nn
\end{eqnarray}
for any $\al_3, \al_4 > 0$.  We estimate the sum of the first two terms in (\ref{BO3}) from below similarly.  Choosing now $\al = 2G(2) / G(1)$ in (\ref{BQ3}) we obtain the lower bound,
\begin{eqnarray} \label{BS3}
&-&  G(2) \bigg\{ 2\al_5 \left< \left( \frac 1{V^2} - \frac 1 4 \right) [(-\De + 2)\Phi]^2 \right> \\
&+& \al_6 \left< (\De\Phi)^2 \right> + \left[ \frac 2{\al_5} + \frac 1{\al_6} + \frac{G(2)}{G(1)} \right] \left< (\na \Phi)^2 \right>\bigg\} \nn \\ 
&-& \frac {G(1)} 8 \left< (\tau_1 \; V-V)^2 \Phi^2 \right> - \frac {G(-1)} 8 \left< (\tau_{-1} \; V-V)^2 \Phi^2 \right> , \nn
\end{eqnarray}
for any $\al_5, \al_6 > 0$.

We may now obtain a lower bound for (\ref{BA3}) by adding (\ref{BK3}) to the final term in (\ref{BO3}) and the expressions of (\ref{BR3}) and (\ref{BS3}).  We obtain the lower bound,
\begin{eqnarray} \label{BT3}
&\ & \left< \left( \frac 1{V^2} - \frac 1 4 \right) \left[ (-\De + 2)\Phi\right]^2 \right> \bigg\{ \frac 3 2 + \frac 1 2 \; G(1) - G(0) \\
&-& 2\al_1 G(1) - 2\al_5 G(2) - 2\al_3 [1 - G(0) - 2G(1)] \bigg\} \nn \\
&+& \frac 1 4 \left< (\De \Phi)^2 \right> \bigg\{ \frac 3 2 + \frac 1 2 \; G(1) - G(0)  - 4\al_2 G(1) - 4\al_6 G(2) \nn \\
&-& 4\al_4 [1 - G(0) - 2G(1)] + 2 [1 - G(0) - 2G(1)] \bigg\} \nn \\
&+& \left< (\na \Phi)^2 \right>  \bigg\{ \frac 3 2 - 2G(1) - \frac 2{\al_1} G(1) - \frac 1{\al_2} G(1) - G(2)^2/G(1) \nn \\
&-& \left[ \frac 2{\al_5} + \frac 1 {\al_6}\right] G(2) - \left[ \frac 2{\al_3} + \frac 1{\al_4} + 2\right] [1-G(0)-2G(1)] \nn \\
&+& [1 - G(0) - 2G(1)] \bigg\} + \frac 1 4 G(0) \left< (4-V^2) \Phi^2 \right> . \nn
\end{eqnarray}

We may rewrite the coefficient of the first term in (\ref{BT3}) as 
\begin{eqnarray} \label{BU3}
\frac 3 2 &+& \frac 1 2 \; G(1) - G(0)  - 2\al_1 G(1) - 2\al_5 G(2) -2\al_3 [1 - G(0) - 2G(1)] \\
= \frac 1 2 &+& (1-2\al_3) [1 - G(0) - 2G(1)] + \left[ \frac 5 2 - 2\al_1 - 2\al_5 \frac{G(2)}{G(1)} \right] G(1). \nn
\end{eqnarray}
If we set now
\begin{equation} \label{BV3}
\al_3 = 1/2, \ \ \al_1  + \al_5 G(2) /G(1) = 5/4,
\end{equation}
we see as in (\ref{BL3}) that the coefficient of the first term in (\ref{BT3}) is $1/2$.  We similarly rewrite the coefficient of the second term as
\begin{equation} \label{BW3}
\frac 3 2 +\frac 1 2 \; G(1) - G(0)  - 4\al_2 G(1) - 4\al_6 G(2)  
\end{equation}
$$
- 4\al_4 [1 - G(0) - 2G(1)] + 2 [1 - G(0) - 2G(1)]  
$$
$$
=\frac 1 2 +(3-4\al_4 )[1 - G(0) - 2G(1)] + \left[ \frac 5 2 - 4\al_2 - 4\al_6 \frac{G(2)}{G(1)} \right] G(1). 
$$
Hence if we set
\begin{equation} \label{BX3}
 \al_2 = 5/8, \ \  \al_4 = 3/4, \ \ \al_6 = 1/8 \; G(2)
\end{equation}
then the second term in (\ref{BT3}) is zero.  We consider the third term in (\ref{BT3}).  This can be written as
\begin{eqnarray} \label{BY3}
&\ & \frac 3 2 - 2G(1) - \frac 2{\al_1} G(1) - \frac 1{\al_2} G(1) - G(2)^2/G(1)  \\
&-& \left[ \frac 2{\al_5} + \frac 1{\al_6}\right] G(2) - \left[ \frac 2{\al_3} + \frac 1{\al_4} + 1\right] [1-G(0)-2G(1)] \nn \\
&=& - \frac 3 2(-\De + 2)G(1) + \left[ \frac 1 2 - \frac 2{\al_3} - \frac 1{\al_4}\right][1 - G(0) - 2G(1)] \nn \\
&+& \left[ 7 - \frac 2{\al_1} - \frac 1{\al_2}\right]G(1)-G(2)^2/G(1) - \left[ \frac 3 2 + \frac 2{\al_5} + \frac 1{\al_6}\right] G(2). \nn
\end{eqnarray}
Using the inequalities (\ref{CB3}) we see this is bounded below by the expression,
\begin{equation} \label{BZ3}
\left\{ 7 - \frac 2{\al_1} - \frac 1{\al_2} + \frac 1 2 \left[ \frac 1 2 - \frac 2{\al_3} - \frac 1{\al_4}\right] - \frac 1{25} - \frac 1 5
\left[ \frac 3 2 + \frac 2{\al_5} + \frac 1{\al_6}\right] \right\}G(1)
\end{equation}
\[   = \left\{ 6.91 - \frac 2{\al_1} - \frac 1{\al_2} - \frac 1{\al_3} - \frac 1{2\al_4}  - \frac 2{5\al_5} - \frac 1{5\al_6} \right\}G(1).
\]
If we substitute the values  (\ref{BV3}), (\ref{BX3}) for $\al_3, \al_4, \al_2, \al_6$ into (\ref{BT3}) we see that the coefficient of $G(1)$ is bounded below by
\begin{equation} \label{CA3}
2.64 - 1.6\; G(2) - g_a(\al_1), \ \ \ a = G(2)/G(1),
\end{equation}
where the function $g_a(z)$ is defined by
\[
g_a(z) = \frac 2 z + \frac{8a}{5[5 - 4z]}, \ \ 0 < z < 5/4, \ \ a > 0.
\]
We can easily compute the minimum of $g_a$ to be
\[	\inf_{0<z< 5/4} g_a(z) = 8\left[ \sqrt{5} + \sqrt{a} \right]^2 / 25.	\]
Using the fact that $a < 1/5, \ G(2) < G(1)/5 < 1/10$ we see that the quantity (\ref{CA3}) is bounded below by $2.64 - .16 - 2.304  > 0$.
\end{proof}
\begin{proof}[Proof of Theorem 1.2]  $(d=2, \ L_1 = 4):$ We need only verify that the function $G$ defined by (\ref{AF3}), (\ref{AZ3}) satisfies the inequalities (\ref{CB3}).  Since $G(y), \ y \ge 1$, decays exponentially one can verify these inequalities with aid of a computer.  In particular we see that
\[	G(0) = .7071, \ \ \ G(1) = .1213, \ \ \ G(2) = .0208,  \]
correct to 4 decimal places, whence (\ref{CB3}) holds. 
\end{proof} 

\section{Proof of Theorem 1.2}
In this section we obtain the formula (\ref{G1}) of $\S$1 for the effective diffusion constant which generalizes the formulas obtained in $\S$3.  We take $L_1 = 2L$ with $L \ge 2$ in Lemma 2.6.  Then $\hat \Omega = \{ (n,y) : 1 \le n \le L, \ y \in \Om_{d-1}\}$.  For $y \in \Om_{d-1}, \ 1 \le j \le L$ we define $\de_j(y), \bar\de_j(y)$ by
\begin{equation} \label{A4}
\de_j(y) = \frac 1{2d} - b(j,y), \ \bar\de_j(y)=\frac 1{2d} + b(j,y).
\end{equation}
We see from (\ref{BK2}), (\ref{A4}) that $\vp^*$ satisfies the system of equations,
\begin{equation}  \label{B4}
\begin{array} {lcccccc}  
\frac{(-\De + 2)}{2d} \vp^*(1,y) &-& \de_2(y) \vp^*(2,y)& -& \de_1(y)\vp^*(1,y) &=& 0, \\
\frac{(-\De + 2)}{2d} \vp^*(2,y)& -& \de_3(y) \vp^*(3,y) &-& \bar\de_1(y)\vp^*(1,y) &=& 0, \\
\cd \cd \cd \cd \cd \cd\cd \cd \cd \cd \cd& -& \cd \cd \cd\cd \cd \cd  & -& \cd \cd \cd\cd \cd \cd &= & 0, \\
\frac{(-\De + 2)}{2d} \vp^*(L-1,y) &-& \de_L(y) \vp^*(L,y) &-& \bar\de_{L-2}(y)\vp^*(L-2,y) &=& 0,  \\
\frac{(-\De + 2)}{2d} \vp^*(L,y)& -& \bar\de_L(y) \vp^*(L,y)& -& \bar\de_{L-1}(y)\vp^*(L-1,y) &=& 0, 
\end{array}
\end{equation}
where $\De = \De_{d-1}$ is the $d-1$ dimensional Laplacian.  If we add all the equations in (\ref{B4}) we obtain the equation
\[  -\De \sum^L_{j=1} \vp^*(j,y) = 0, \ \ y \in \Om_{d-1}.		\]
On using the normalization $\left<\vp^*\right>_{\hat \Om} = 1$ we conclude that
\begin{equation} \label{C4}
\sum^L_{j=1} \vp^*(j,y) = L, \ \ y \in \Om_{d-1}.
\end{equation}
Evidently we can rewrite the first equation of (\ref{B4}) as
\begin{equation} \label{D4} 
\left[ -\frac{\De}{2d} + \bar \de_1(y) \right] \vp^*(1,y) - \de_2(y)\vp^*(2,y) = 0.
\end{equation}
If we add (\ref{D4}) to the second equation of (\ref{B4}) we obtain the equation
\begin{equation} \label{E4}
\left( -\frac{\De}{2d}\right)\vp^*(1,y) +  \left[- \frac{\De}{2d} + \bar \de_2(y) \right] \vp^*(2,y) - \de_3(y)\vp^*(3,y) = 0.
\end{equation}
Adding (\ref{E4}) to the third equation of (\ref{B4}) and proceeding similarly with subsequent equations we obtain the system,
\begin{equation} \label{F4}
\left(- \frac{\De}{2d}\right)\vp^*(1,y) + \left(- \frac{\De}{2d}\right)\vp^*(2,y) +  \left[ -\frac{\De}{2d} + \bar \de_3(y) \right] \vp^*(3,y) - \de_4(y)\vp^*(4,y) = 0,
\end{equation}
\centerline{ $\cd \ \ \ \cd \ \ \ \cd \ \ \ \cd \ \ \ \cd \ \ \ \cd \ \ \ \cd \ \ \ \cd \ \ \ \cd \ \ \ $}
\[
\left( -\frac{\De}{2d}\right)\vp^*(1,y) + \cdots + \left(- \frac{\De}{2d}\right)\vp^*(L-2,y) +  \left[- \frac{\De}{2d} + \bar \de_{L-1}(y) \right] \vp^*(L-1,y) - \de_L(y)\vp^*(L,y) = 0, \]
where we have omitted the final equation of (\ref{B4}).  From (\ref{D4}), (\ref{E4}), (\ref{F4}) we can write $\vp^*(j,y), \ 2 \le j \le L$, in terms of $\vp^*(1,y)$.  Substituting these into (\ref{C4}) we obtain an equation for $\vp^*(1,y)$ of the form
\begin{equation} \label{G4}
\mathcal{L} \; \vp^*(1,y) = L, \ \ \ y \in \Om_{d-1},
\end{equation}
where $\mathcal{L}$ is an operator on functions on $\Om_{d-1}$.

Next we consider the equations (\ref{AW2}), (\ref{BB2}) for the function $\psi_0$ on $\hat \Om$.  Thus $\psi_0$ satisfies the system of equations, 
\begin{equation} \label{H4}
\begin{array}{lcccccc}
\frac{(-\De + 2)}{2d} \psi_0(1,y) &-& \bar\de_1(y) \psi_0(2,y) &+&\de_1(y)\psi_0(1,y) &=& 0, \\
\frac{(-\De + 2)}{2d} \psi_0(2,y) &-& \bar\de_2(y) \psi_0(3,y)& -& \de_2(y)\psi_0(1,y)& = &0, \\
\cd \cd \cd \cd \cd \cd\cd \cd \cd \cd \cd& -& \cd \cd \cd\cd \cd \cd  & -& \cd \cd \cd\cd \cd \cd &= & 0, \\
\frac{(-\De + 2)}{2d} \psi_0(L-1,y) &-& \bar\de_{L-1}(y) \psi_0(L,y)& -& \de_{L-1}(y)\psi_0(L-2,y) &= &0,  \\
\frac{(-\De + 2)}{2d} \psi_0(L,y) &+& \bar\de_L(y) \psi_0(L,y)& -& \de_L(y)\psi_0(L-1,y) &= &\bar\de_L(y). 
\end{array}
\end{equation}
We can rewrite the first equation of (\ref{H4}) as
\begin{equation} \label{J4}
\frac{(-\De + 4)}{2d} \psi_0(1,y) = \bar\de_1(y) \left[ \psi_0(1,y) + \psi_0(2,y)\right] = \bar\de_1(y) u(1,y),  \  \ y \in \Om_{d-1},
\end{equation}
where $u(1) = \psi_0(1) + \psi_0(2)$.  If we put now $u(2) = \psi_0(3) - \psi_0(1)$ then on using $(\ref{J4})$ we see that the second equation of (\ref{H4}) is the same as
\begin{equation} \label{K4}
\left[ -\frac{\De}{2d} + \de_1\right]u(1) - \bar \de_2u(2) = 0.
\end{equation}
Observe that (\ref{K4}) is identical to (\ref{D4}) under the reflection $b \ra -b$.  We can similarly obtain the reflection of the equations (\ref{E4}), (\ref{F4}) by defining the variables $u(j), j = 3,...$ by
\[	u(j) = \psi_0(j+1) - \psi_0(j-1),  \  \ j = 3,..., L-1.	\]
Let us assume that the $u(j), j=1,.....,J-1$, satisfy the reflection of the first $J-2$ of the equations (\ref{D4}), (\ref{E4}), (\ref{F4}).   We show that $u(J)$ then satisfies the $(J-1)$st equation provided $J \le L-1$.  To see this we consider the $J$th equation of (\ref{H4}) which we may write  as 
\begin{equation} \label{L4}
\left( \frac{-\De + 2}{2d}\right) \psi_0(J) - \bar\de_J \left[ u(J) + \psi_0 (J-1) \right] - \de_J \psi_0(J-1) = 0.
\end{equation}
We may rewrite (\ref{L4}) as
\begin{multline} \label{M4}
\left[- \frac{\De}{2d} + \de_{J-1}\right] u(J-1)  - \bar\de_J  u(J) +\bar \de_{J-1} u(J-1) + \\
 \frac{(-\De + 2)}{2d} \psi_0(J-2) - \frac 2{2d} \psi_0(J-1)=0.
\end{multline}

If $J=3$ then we have that
\begin{eqnarray*}
\bar \de_{J-1} u(J-1) &+& \frac{(-\De + 2)}{2d} \psi_0(J-2) - \frac 2{2d} \psi_0(J-1) = \\
\bar \de_2 u(2) &+& \frac{(-\De + 4)}{2d} \psi_0(1) - \frac 2{2d} \left[ \psi_0(1) + \psi_0(2) \right] = \\
\bar \de_2 u(2) &+& \bar \de_1 u(1) - \frac 2{2d} u(1) = \bar \de_2 u(2) - \de_1 u(1) \\
&=& \left[ -\frac{\De}{2d} + \de_1\right] u(1) - \de_1u(1) =  -\frac{\De}{2d} u(1).
\end{eqnarray*}
We have shown that the result holds for $J=3$.  More generally we have that
\begin{eqnarray*}
\bar \de_{J-1} u(J-1) &+& \frac{(-\De + 2)}{2d} \psi_0(J-2) - \frac 2{2d} \psi_0(J-1) = \left[- \frac{\De}{2d} + \de_{J-2}\right] u(J-2)\\
&+& \sum^{J-3}_{j=1} - \frac{\De}{2d}  u(j) + \frac{(-\De + 2)}{2d} \psi_0(J-2) - \frac 2{2d} \psi_0(J-1) \\
&=& \sum^{J-2}_{j=1} - \frac{\De}{2d} u(j) - \bar\de_{J-2} u(J-2) - \frac 2{2d}\psi_0(J-3) \\
&+& \frac{(-\De + 2)}{2d} \psi_0 (J-2).
\end{eqnarray*}

To complete the proof we need then to show that
\[  -\bar\de_{J-2} \; u(J-2) - \frac2{2d} \psi_0(J-3) + \frac{(-\De + 2)}{2d} \psi_0(J-2) = 0.	\]
This last equation is however simply the $(J-2)$nd equation of (\ref{H4}).

We have shown that $u(j), j=1,...,L-1$ satisfies the reflection of the first $L-2$ equations of  (\ref{D4}), (\ref{E4}), (\ref{F4}).  Define now $u(L) = 1 - \psi_0(L) - \psi_0(L-1)$, whence there is the identity,
\begin{equation} \label{N4}
\sum^L_{j=1} \ u(j) = 1.
\end{equation}
 We shall show that the $u(j), j=1,...,L$ satisfy the reflection of the final equation of (\ref{F4}).  To see this we write the final equation of (\ref{H4}) as
\[ \left( \frac{-\De + 2}{2d}\right) \psi_0(L) + \bar\de_L \left[ 1-u(L) - \psi_0 (L-1) \right] - \de_L \psi_0(L-1) = \bar\de_L, \]
whence we have that
\[ \left( \frac{-\De + 2}{2d}\right) \left[ u(L-1) + \psi_0 (L-2) \right] - \bar \de_L u(L) - \frac 2{2d} \psi_0(L-1) = 0. \]
We may rewrite the previous equation as
\[
\left[ -\frac{\De}{2d} + \de_{L-1}\right] u(L-1) - \bar \de_L \; u(L) + \bar\de_{L-1}  u(L-1) + \frac{(-\De + 2)}{2d} \psi_0 (L-2) - \frac 2{2d} \psi_0(L-1) = 0.    \]

Now if we use the identity already established,
\[     \bar \de_{L-1} u(L-1) =  \left[ -\frac{\De}{2d} + \de_{L-2}\right] u(L-2) + \sum^{L-3}_{j=1}  - \frac{\De}{2d} u(j),  \]
we see that it is sufficient to show that
\[  \de_{L-2} u(L-2) + \frac{(-\De + 2)}{2d} \psi_0(L-2) - \frac 2{2d} \psi_0(L-1) = 0.    \]
This last equation is just the $(L-2)$nd equation of (\ref{H4}).

Let $\mathcal{L}_R$ be the reflection of the operator $\mathcal{L}$ of (\ref{G4}) obtained by replacing $b$ by $-b$.  Then on comparing (\ref{C4}), (\ref{N4}) we see that $u(1)$ satisfies the equation
\begin{equation} \label{O4}
\mathcal{L}_R \ u(1) = 1.
\end{equation}
We are able now to come up with a new formula for the effective diffusion constant.  On using (\ref{BC2}), (\ref{G4}), (\ref{J4}), (\ref{O4}) we have that the effective diffusion constant is given by
\begin{equation} \label{P4}
8 L^2 d \left< \left[\de_1 \mathcal{L}^{-1} \; 1\right] (-\De + 4)^{-1} \left[\bar\de_1 \mathcal{L}^{-1}_R \; 1\right] \right>,
\end{equation}
where $\left< \; \cd \; \right>$ is the uniform probability measure on $\Om_{d-1}$.  The formula (\ref{G1}) follows from (\ref{P4}). In order for (\ref{P4}) to be valid we need to show that $\mathcal{L}$ is invertible.
\begin{lemma} Let $\mathcal{L}$ be the matrix defined by (\ref{G4}).  Then $\mathcal{L}$ is invertible and the matrix $\mathcal{L}^{-1}$ has all positive entries.
\end{lemma}
\begin{proof} We proceed by induction.  For $k = 2,3.....$ let $\mathcal{L}_k$ be the operator (\ref{G4}) when $L = k$.  It is easy to see from (\ref{C4}) - (\ref{F4}) that the $\mathcal{L}_k$ satisfy the recurrence relation,
\begin{equation} \label{Q4}
\mathcal{L}_{k+1} = \frac 1{\de_{k+1}} \left[ -\frac{\De}{2d} + \bar \de_k + \de_{k+1} \right] \mathcal{L}_k - \frac{\bar\de_k}{\de_{k+1}} \mathcal{L}_{k-1}, \  k \ge 1; \ \mathcal{L}_0 = 0, \  \mathcal{L}_1 = 1.
\end{equation}
The result will follow by showing that the matrices $A_k = \mathcal{L}_{k-1} \mathcal{L}_k^{-1}, k \ge 2$, have all positive entries and principal eigenvalue strictly less than 1. Evidently this is the case for $k=2$. Now from (\ref{Q4}) we see that the $A_k$ satisfy the recurrence relation,
\begin{equation} \label{R4}
A_{k+1} = \left\{ -\frac{\De}{2d} + \bar \de_k + \de_{k+1} - \bar \de_k A_k\right\}^{-1} \de_{k+1}.
\end{equation} 
If $A_k$ has all positive entries with principal eigenvalue strictly less than 1 then the matrix $[-\De/2d + \bar \de_k + \de_{k+1}]^{-1}\bar \de_k A_k$ has the same property and the matrix $A_{k+1}$ defined by (\ref{R4}) has all positive entries.  To conclude the induction step we need therefore to show that $A_{k+1}$  has principal eigenvalue strictly less than $1$.  To see this note that if $1$ denotes the vector with all entries $1$ then  
\[     \left\{- \frac{\De}{2d} + \bar \de_k + \de_{k+1} - \bar \de_k A_k\right\} 1 > \de_{k+1},    \]
whence we conclude that
\[    \left\{ -\frac{\De}{2d} + \bar \de_k + \de_{k+1} - \bar \de_k A_k\right\}^{-1} \de_{k+1} (1) < 1.   \]
\end{proof}

\thanks{ {\bf Acknowledgement:} This research was 
partially supported by NSF under grant DMS-0500608.

\end{document}